\documentclass[12pt, leqno]{amsart}

\usepackage[OT2,T1]{fontenc}
\usepackage{indentfirst}
\usepackage{amstext}
\usepackage{amsthm}
\usepackage{amsopn}
\usepackage{amsfonts}
\usepackage{amsmath}
\usepackage{latexsym}
\usepackage[francais,english]{babel}
\usepackage{amscd}
\usepackage{amssymb}
\usepackage{amsmath}
\usepackage[all,cmtip]{xy}
\usepackage{leftidx}
\usepackage{graphicx}
\usepackage{etoolbox}
\patchcmd{\section}{\normalfont\scshape\centering}{\normalfont\bfseries}{}{}
\patchcmd{\subsection}{-.5em}{.5em}{}{}

\newtheorem{theo}{{Theorem}}[section]
\newtheorem{coro}[theo]{{Corollary}}
\newtheorem{lemma}[theo]{{Lemma}}
\newtheorem{prop}[theo]{Proposition}

\theoremstyle{definition}
\newtheorem{remark}[theo]{\textbf{Remark}}
\newtheorem{defn}[theo]{Definition}

\numberwithin{equation}{section}

\newcommand{\ra}{\rightarrow}

\newcommand{\gm}{\mathbb{G}}

\DeclareFontEncoding{OT2}{}{} % to enable usage of cyrillic fonts
  \newcommand{\textcyr}[1]{%
    {\fontencoding{OT2}\fontfamily{wncyr}\fontseries{m}\fontshape{n}%
     \selectfont #1}}
\newcommand{\sha}{{\mbox{\textcyr{Sh}}}}

\def\YEAR{\year}\newcount\VOL\VOL=\YEAR\advance\VOL by-1995
\def\firstpage{1}\def\lastpage{1000}
\def\received{}\def\revised{}
\def\communicated{}

\makeatletter
\def\magnification{\afterassignment\m@g\count@}
\def\m@g{\mag=\count@\hsize6.5truein\vsize8.9truein\dimen\footins8truein}
\makeatother

\font\eightrm=cmr8
\font\caps=cmcsc10                    % Theorem, Lemma etc
\font\Caps=cmcsc10 scaled \magstep1   % Title
\font\scaps=cmcsc8

\def\bfseries{\normalsize\caps}

\makeatletter
\@specialpagefalse\headheight=8.5pt
\def\DocMath{{\def\th{\thinspace}\scaps Documenta Math.}}
\renewcommand{\@oddfoot}{\hfill\scaps Documenta Mathematica 
    \number\VOL\  (\number\YEAR) \number\firstpage--\lastpage\hfill}
\renewcommand{\@evenfoot}{\ifnum\thepage>\lastpage\hfill\scaps
    Documenta Mathematica \number\VOL\  (\number\YEAR)\hfill\else\@oddfoot\fi}%

\renewcommand{\@evenhead}{%
    \ifnum\thepage>\lastpage\rlap{\thepage}\hfill%
    \else\rlap{\thepage}\slshape\leftmark\hfill{\caps\SAuthor}\hfill\fi}%
\renewcommand{\@oddhead}{%
    \ifnum\thepage=\firstpage{\DocMath\hfill\llap{\thepage}}%
    \else{\slshape\rightmark}\hfill{\caps\STitle}\hfill\llap{\thepage}\fi}%
\makeatother

\def\TSkip{\bigskip}
\newbox\TheTitle{\obeylines\gdef\GetTitle #1
\ShortTitle  #2
\SubTitle    #3
\Author      #4
\ShortAuthor #5
\EndTitle
{\setbox\TheTitle=\vbox{\baselineskip=20pt\let\par=\cr\obeylines%
\halign{\centerline{\Caps##}\cr\noalign{\medskip}\cr#1\cr}}%
	\copy\TheTitle\TSkip\TSkip%
\def\next{#2}\ifx\next\empty\gdef\STitle{#1}\else\gdef\STitle{#2}\fi%
\def\next{#3}\ifx\next\empty%
    \else\setbox\TheTitle=\vbox{\baselineskip=20pt\let\par=\cr\obeylines%
    \halign{\centerline{\caps##} #3\cr}}\copy\TheTitle\TSkip\TSkip\fi%
\centerline{\caps #4}\TSkip\TSkip%
\def\next{#5}\ifx\next\empty\gdef\SAuthor{#4}\else\gdef\SAuthor{#5}\fi%
\ifx\received\empty\relax
    \else\centerline{\eightrm Received: \received}\fi%
\ifx\revised\empty\TSkip%
    \else\centerline{\eightrm Revised: \revised}\TSkip\fi%
\ifx\communicated\empty\relax
    \else\centerline{\eightrm Communicated by \communicated}\fi\TSkip\TSkip%
\catcode'015=5}}\def\Title{\obeylines\GetTitle}
\def\Abstract{\begingroup\narrower
    \parskip=\medskipamount\parindent=0pt{\caps Abstract. }}
\def\EndAbstract{\par\endgroup\TSkip}

\long\def\MSC#1\EndMSC{\def\arg{#1}\ifx\arg\empty\relax\else
     {\par\narrower\noindent%
     2010 Mathematics Subject Classification: #1\par}\fi}

\long\def\KEY#1\EndKEY{\def\arg{#1}\ifx\arg\empty\relax\else
	{\par\narrower\noindent Keywords and Phrases: #1\par}\fi\TSkip}

\newbox\TheAdd\def\Addresses{\vfill\copy\TheAdd\vfill
    \ifodd\number\lastpage\vfill\eject\phantom{.}\vfill\eject\fi}
{\obeylines\gdef\GetAddress #1
\Address #2 
\Address #3
\Address #4
\EndAddress
{\def\xs{4.3truecm}\parindent=0pt
\setbox0=\vtop{{\obeylines\hsize=\xs#1\par}}\def\next{#2}
\ifx\next\empty % 1 address
     \setbox\TheAdd=\hbox to\hsize{\hfill\copy0\hfill}
\else\setbox1=\vtop{{\obeylines\hsize=\xs#2\par}}\def\next{#3}
\ifx\next\empty % 2 addresses
     \setbox\TheAdd=\hbox to\hsize{\hfill\copy0\hfill\copy1\hfill}
\else\setbox2=\vtop{{\obeylines\hsize=\xs#3\par}}\def\next{#4}
\ifx\next\empty\ % 3 addresses
     \setbox\TheAdd=\vtop{\hbox to\hsize{\hfill\copy0\hfill\copy1\hfill}
                \vskip20pt\hbox to\hsize{\hfill\copy2\hfill}}
\else\setbox3=\vtop{{\obeylines\hsize=\xs#4\par}}
     \setbox\TheAdd=\vtop{\hbox to\hsize{\hfill\copy0\hfill\copy1\hfill}
	        \vskip20pt\hbox to\hsize{\hfill\copy2\hfill\copy3\hfill}}
\fi\fi\fi\catcode'015=5}}\gdef\Address{\obeylines\GetAddress}

\def\LOCAL{\jobname.files}

\begin{document}
%%%%% ------------- fill in your data below this line  -------------------
%%%%%    The following lines \Title ... \EndAddress must ALL be present
%%%%%    and in the given order.
\Title {On unramified Brauer groups of  torsors over tori}
%%%%%    Put here the title. Line breaks will be recognized. 
\ShortTitle  
%%%%%    Running title for odd numbered pages, ONE line, please. 
%%%%%    If none is given, \Title will be used instead.          
\SubTitle   
%%%%%    A possible subtitle goes here.
\Author E. Bayer-Fluckiger and  R. Parimala
%%%%%    Put here name(s) of authors. Line breaks will be recognized.  
\ShortAuthor 
%%%%%%   Running title for even numbered pages, ONE line, please. 
%%%%%%   If none is given, \Author will be used instead.          
\EndTitle
\Abstract 
%%%%%    Put here the abstract of your manuscript.
In this paper we introduce a method to obtain algebraic information using arithmetic one in 
the study of tori and their principal homogeneous spaces. In particular, using some results of 
he authors with Ting-Yu Lee, we determine the unramified Brauer groups of some norm one tori, and their torsors. 

\EndAbstract
\MSC  11G35, 12G05,14G05, 14G25.
%%%%%    2010 Mathematics Subject Classification: 

\EndMSC
\KEY 
%%%%%    Keywords and Phrases:     
Unramified Brauer groups, Tate-Shafarevich groups, Hasse principle and weak approximation, norm tori.
\EndKEY
%%%%%    All 4 \Address lines below must be present. To center the last
%%%%%    entry, no empty lines must be between the following \Address
%%%%%    and \EndAddress lines.
\Address 
%%%%%    Address of first Author here
EPFL-FSB-MATH\\ %
Station 8 \\ %
1015 Lausanne \\ %
Switzerland \\ %
eva.bayer@epfl.ch 

\Address
%%%%%    Address of second Author here etc.
Department of Mathematics  \\ %
Emory University \\ %
400 Dowman Drive~NE \\ %
Atlanta, GA 30322, USA \\ %
praman@emory.edu

\Address
\Address
\EndAddress
%%
%%       Make sure the last tex command in your manuscript
%%       before the first \end{document} is the command  \Addresses
%%
%%---------------------Here the prologue ends---------------------------------
%%--------------------Here the manuscript starts------------------------------

%%--------------------Here the manuscript ends--------------------------------
\Addresses

\tolerance 400 \pretolerance 200 \selectlanguage{english}

%\title{On unramified Brauer groups of  torsors over tori}
%\author{E. Bayer-Fluckiger}
%\address{EPFL-FSB-MATH\\ %
%Station 8 \\ %
%1015 Lausanne \\ %
%Switzerland}
%\email{eva.bayer@epfl.ch}
%
%
%\author{R. Parimala}
%\address{Department of Mathematics  \\ %
%Emory University \\ %
%400 Dowman Drive~NE \\ %
%Atlanta, GA 30322, USA}
%\email{praman@emory.edu}
%
%
%%\date{\today}
%\maketitle
%
%\begin{abstract}
%In this paper we introduce a method to obtain algebraic information using arithmetic one in the study of tori and their principal homogeneous spaces. In particular, using some results of the authors with Ting-Yu Lee, we determine the unramified Brauer groups of some norm one tori, and their torsors. 
%
%\medskip
%
%\noindent {\em Keywords:} Unramified Brauer groups, Tate-Shafarevich groups, Hasse principle and weak approximation, norm tori.
%
%
%
%
%\noindent {\em MSC 2000:} 11G35, 12G05,14G05, 14G25.
%\end{abstract}
%

\small{} \normalsize

\medskip

\selectlanguage{english}
\section{Introduction}
In this paper, we use arithmetic information to obtain algebraic ones. Let $G$ be a finite group, and let $M$ be a $G$-lattice. Let $\ell'/\ell$ a finite
{\it unramified} extension of number fields with Galois group $G$; such an extension exists by \cite{F}.
Let $T$ be an $\ell$-torus with character group $M$. We have the following isomorphisms
$$(*) \ \ \ \ \sha^2_{\rm cycl}(G,M) \simeq \sha^2(\ell,M) \simeq \sha^1(\ell,T)^*;$$
the first isomorphism is Proposition \ref{lemma}, the second one follows from Poitou-Tate duality (see \S \ref {alg} and \S \ref{arith} for the notation).  In the following, we combine  $(*)$ with arithmetic results as well as
some theorems of Colliot-Th\'el\`ene and Sansuc; we illustrate the results with the following example (see \S \ref {un cyclic}) :

\medskip
\noindent
{\bf Example.}
Let $k$ be a field, and let
$L = K_1 \times \dots \times K_n$, where $K_1, \dots, K_n$ are cyclic extensions of $k$ of
prime degree $p$.   Let  ${\rm N}_{L/{k}}: L \to k$ denote the norm map, and let
$T_{L/{k}}  = R^{(1)}_{L/k}({\bf G}_m)$ be the $k$-torus
defined by $$1 \to T_{L/{k}} \to R_{L/k}({\bf G}_m) {\buildrel {\rm N}_{L/{k}} \over \longrightarrow}  {\bf G}_m \to 1.$$ 
Let $k'/k$ be a Galois extension of minimal degree splitting $T_{L/{k}}$, and let $G = {\rm Gal}(k'/k)$. Set $T = T_{L/{k}}$, let $T^c$ be a smooth
compactification of $T$, and let ${\rm Br}(T^c)$ be its Brauer group.
%Let $M$ be the $G$-lattice of characters of $T_{L/{k}}$. 

%\medskip Set  $I = \{1, \dots, n\}$. 
%We consider the norm polynomials ${\rm N}_{K_i/k}(t_i)$  for $i \in I$ as elements of $k[X]$.  For all $i \in I$, set $N_i = {\rm N}_{K_i/{k}}(t_i)$ and let $\sigma_i$ be a generator
%of ${\rm Gal}(K_i/k)$.  We denote by $(N_i,K_j)$ the cyclic algebra over $k(X^c)$ associated
%to the generator $\sigma_j$. 

\medskip
Let $a \in k^{\times}$, and let $X$ be the affine $k$-variety determined by the equation ${\rm N}_{L/k}(t) = a$; $X$ is a torsor under $T_{L/k}$.
Let $X^c$ be
a smooth compactification of $X$. We denote by ${\rm Br}(X^c)$ the Brauer group of $X^c$.

\medskip
\noindent
{\bf Theorem.} {\it 
{\rm (a)} If $G \not \simeq C_p \times C_p$, then $${\rm Br}(T^c)  / {\rm Br}(k) = {\rm Br}(X^c)/{\rm Im}({\rm Br}(k)) = 0.$$

%\smallskip
%{\rm (b)} If $G  \simeq C_p \times C_p$, then $${\rm Br}(T^c)  / {\rm Br}(k) \simeq ({\bf Z} / p {\bf Z})^{n-2}.$$

\medskip
{\rm (b)} If $G  \simeq C_p \times C_p$, then $${\rm Br}(T^c)  / {\rm Br}(k) \simeq {\rm Br}(X^c)/{\rm Im}({\rm Br}(k))  \simeq ({\bf Z} / p {\bf Z})^{n-2}.$$} 

%\medskip
%{\rm (c)} If $G  \simeq C_p \times C_p$, then ${\rm Br}(X^c)/{\rm Br}(k) $ is generated by the classes of the elements $$(N_i,K_n), \ \ i = 1,\dots,n-2.$$ in ${\rm Br}(k(X^c))$. }

\medskip This is proved in Theorem \ref{brun} %and \ref{brun p}, 
using $(*)$ as well as some (arithmetic) results of \cite{BLP}. 

%\medskip
%Part (c) is proved by a combination of arithmetic and algebraic methods. Indeed,  ${\rm Br}(X^c)  / {\rm Im}({\rm Br}(k) )$ injects
%into ${\rm Br}(T^c)/{\rm Br}(k)$ (see \S \ref {un}). We then use algebraic methods to show that  if ${\rm char}(k) \not = p$, then ${\rm Br}(X^c)  / {\rm Im}({\rm Br}(k) )$
%contains $n - 2$ elements that are linearly independent over ${\bf Z} / p {\bf Z})$, hence we obtain part (c) of the above theorem (see Theorems \ref{brun p}
%and Theorem \ref{generators}).

\medskip Further, we also give generators for the group ${\rm Br}(X^c)/{\rm Im}({\rm Br}(k))$ (see Theorem \ref{generators}), in the spirit of Colliot-Th\'el\`ene's results
for biquadratic extensions (see \cite {CT}, \S 4). 

\medskip
The starting point is the following key observation of Jean-Louis Colliot-Th\'el\`ene :

\medskip
\noindent
{\bf Proposition.} {\it Let $G$ be a finite group, and let $M$ be a $G$-lattice. If for all number fields $k$ and every $k$-torus $T$ with character group isomorphic to
the Galois module $M$ via a surjection $\Gamma_k \to G$ one can show that $\sha^1(k,T) = 0$, then $\sha_{\rm cycl}^2(G,M) = 0$, and  every such $k$-torus $T$ has weak approximation.}

\medskip  Since every finite group is the Galois group of some unramified extension of number fields, we may
realize the purely algebraic group $\sha_{\rm cycl}^2(G,M)$ as the Tate-Shafarewich group of a torus over a number field; this is summarized in $(*)$.
If the module satisfies the hypothesis of the proposition, then $\sha_{\rm cycl}^2(G,M) = 0$, and weak approximation follows from an exact
sequence due to Voskresenskii (see \ref {Vosk}). As we will see, the hypotheses of the proposition can be weakened (see Corollary
\ref {hp wa}).

\medskip
\noindent
We thank Jean-Louis Colliot-Th\'el\`ene for sharing his ideas with us, as well as for several useful suggestions.\\
\noindent
The second author is partially supported by the National Science Foundation Grant DMS--1801951.
 
%$a \geqslant b$

\section{Algebraic preliminaries}\label{alg}

\medskip
Let $k$ be a field, let $k_s$ be a separable closure of $k$ and set $\Gamma_k = {\rm Gal}(k_s/k)$.
 We fix once and for all this
separable closure $k_s$, and all separable extensions of $k$ that will appear in the paper will be contained in $k_s$.
We use standard notation in Galois cohomology; in particular, if $M$ is
a discrete $\Gamma_k$-module and $i$ is an integer  $\ge 0$, we set $H^i(k,M) = H^i(\Gamma_k,M)$. A {\it $\Gamma_k$-lattice} will
be a torsion free  $\bf Z$-module of finite rank on which $\Gamma_k$ acts continuously.

\medskip
\begin{lemma}\label{1.3} Let $M$ be a $\Gamma_k$-lattice, and let $k'/k$ be a finite Galois extension with Galois group $G$
such that $\Gamma_{k'}$ acts trivially on $M$. Then the natural map $H^2(G,M) \to H^2(k,M)$ has trivial kernel.

\end{lemma}

\noindent
{\bf Proof.} Since $M$ is isomorphic to the trivial $\Gamma_{k'}$-module ${\bf Z}^n$
 for some $n$, we
 have $H^1(k',M) = 0$. Hence the
 exact sequence of groups $0 \to \Gamma_{k'} \to \Gamma_{k} \to G \to 0$ induces an exact sequence in Galois cohomology (cf. \cite{Se}, page 118, Prop. 5)
$$(*) \hskip 0.5 cm 0 \to  H^2(G,M) \to H^2(k,M) \to H^2(k',M)^G.$$ 
%$$(*) \hskip 1cm H^1(k',M) \to H^2(G,M) \to H^2(k,M).$$ 
Therefore the map $H^2(G,M) \to H^2(k,M)$ has trivial kernel, as claimed. 
 
 \medskip
 Let $G$ be a finite group. A {\it $G$-lattice} is by definition a ${\bf Z}$-torsion free ${\bf Z}[G]$-module of finite rank. For a $k$-torus $ T$, we denote by $\hat{ T} = {\rm Hom}(T,\gm_m)$ its character group; it is a $\Gamma_k$-lattice.

 \begin{prop}\label{equivalence of categories} Let $M$ be a $G$-lattice. Let $\eta: \Gamma_k \to G$ be a continuous surjective homomorphism. There exists a $k$-torus $T$ such that 
 $\hat T$ is isomorphic to the $G$-lattice $M$, regarded as $\Gamma_k$-lattice through $\eta$.
 %the surjection $\Gamma_k \to G$. 

\end{prop}

\noindent 
{\bf Proof.} See \cite {Bo}, Chapter III, 8.12.

\medskip
If $g \in G$, we denote by $\langle g \rangle$ the cyclic subgroup of $G$ generated by $g$. Let $M$ be
a $G$-lattice. The cyclic Tate-Shafarevich group  $\sha_{\rm cycl}^2(G,M)$ is the group

$$\sha^2_{cycl}(G,M) = {\rm Ker}[H^2(G,M) \to \prod_{g \in G} H^2(\langle g \rangle,M)].$$

\medskip
We recall a result of Colliot-Th\'el\`ene and Sansuc :

\begin{theo}\label{Colliot Br cycl} Let $G$ be a finite group, let $T$ be a $k$-torus, and assume that the character group of $T$ is a $G$-lattice via a surjection $\Gamma_k \to G$. Let $T^c$ be a smooth compactification of $T$. We have ${\rm Br}(T^c)/{\rm Br}(k) \simeq \sha^2_{\rm cycl}(G,\hat T)$.

\end{theo}

\noindent
{\bf Proof.} See  \cite {CTS 87}, Theorem 9.5 (ii). In \cite {CTS 87} ,  the hypothesis ${\rm char}(k) = 0$ is only used to ensure 
the existence of a a smooth
compactification of $T$; since this is now known in any characteristic (see \cite {CTHSk 05}), the result holds in general.

\medskip

Let $Y$ be a smooth projective, geometrically connected $k$-variety, and set $\overline Y = Y \times_k k_s$.  We have the
following spectral sequence 

$$H^p(k,H^q(\overline Y, {\bf G}_m)) \implies H^{p+q}(Y,{\bf G}_m)$$ giving the exact sequence
$$(**) ~~~~~H^2(k,{\bf G}_m) \to {\rm Ker}[H^2_{\rm et}(Y,{\bf G}_m) \to H^2_{\rm et}(\overline Y,{\bf G}_m)] \to H^1(k,{\rm Pic}(\overline Y) ) \to $$
$$\to {\rm Ker}[H^3_{\rm et}(k,{\bf G}_m) \to H^3_{\rm et}(Y,{\bf G}_m)].$$ 

We refer to \cite{CTHSk 03}, Section 2 for the following theorem.

\begin{theo} 
\label{brnr}
Let  $Y$ be a smooth projective variety defined over $k$ with $\overline{Y}$ birational   to the projective space.
Then there is an injection 
$${\rm Br}(Y)/{\rm Im (Br}(k)) \to  H^1(k,{\rm Pic}(\overline Y)).$$
If further $Y(k) \neq \emptyset$, the above injection yields  an isomorphism  ${\rm Br}(Y)/{\rm Br}(k) \simeq  H^1(k,{\rm Pic}(\overline Y)).$
\end{theo}

\noindent 
{\bf Proof.}
Since $\overline{Y}$ is birational to the projective space, 
 $H^2_{\rm et}(\overline Y,{\bf G}_m)
= {\rm Br}(\overline Y) = 0$ (cf. \cite{CTSk}, Theorem 5.1.3, Corollary 5.2.6) and we have an injection  
$$  {\rm Br}(Y)/{\rm Im (Br}(k)) \to H^1(k,{\rm Pic}(\overline Y)).$$ 
If further $Y(k) \neq  \emptyset$, ${\rm Ker}[H^i (k,{\bf G}_m) \to H^i_{\rm et}(Y, {\bf G}_m)] = 0$ for  all $i$, so that
the injection above becomes
an isomorphism $${\rm Br}(Y)/{\rm Br}(k) \simeq H^1(k,{\rm Pic}(\overline Y)).$$

\section{Arithmetic preliminaries}\label{arith}

\medskip
Let $k$ be a global field, and let
$\Omega_k$ be the set of all places of $k$; if $v \in \Omega_k$, we denote by $k_v$ the completion of $k$ at $v$. 

\medskip
For any $k$-torus $T$, set $\sha ^i(k,T) = {\rm Ker} (H^i(k,T) \ra \underset {v \in \Omega_k} \prod H^i(k_v,T))$. 
If $M$ is a $\Gamma_k$-module,
set 
$\sha^i (k,M) = {\rm Ker} (H^i(k,M) \ra \underset {v \in \Omega_k} \prod H^i(k_v,M))$, and let $\sha_{\omega}^i(k,M)$ be the set of
$x \in H^i(k,M)$ that map to 0 in $H^i(k_v,M)$ for almost all $v \in \Omega_k$. Recall that by  Poitou-Tate duality, we have 
an isomorphism of finite groups
$\sha^2(k,\hat T) \simeq \sha^1(k,T)^{\ast}$, where $\sha^1(k,T)^{\ast}=  {\rm Hom}(\sha^1(k,T), \bf Q / \bf Z)$
 denotes the dual of $\sha^1(k,T)$. We denote by $G_M$ the kernel of the map $\Gamma_k \to {\rm Aut}(M)$; set $k_M = (k_s)^{G_M}$, and
let $G = \Gamma_k/G_M$. If $v \in \Omega_k$, we denote by $G_v$ the decomposition group of  a prime $w$ lying over $v$ in the extension $k_M/k$. For various $w$ over $v$, the groups $G_v$'s are conjugate and  are
subgroups of $G$. Let $\sha^i_{\rm cycl}(k,M)$ be the set of $x \in H^i(k,M)$ that map to 0 in $H^i(k_v,M)$ for all $v \in \Omega_k$ such that 
$G_v$  is cyclic.
\medskip The following lemmas are well-known :

\begin{lemma}\label{Z} Let $\bf Z$ be the trivial $\Gamma_k$-module. Then $\sha^2_\omega(k,{\bf Z}) = 0$.
In particular $\sha^2(k,  {\bf Z}) = 0$.
\end{lemma}

\noindent
{\bf Proof.} Let $L/k$ be a field extension. The trivial $\Gamma_L$-module $\bf Q$ is uniquely divisible, hence $H^i(L,{\bf Q}) = 0$ for all $i \ge 1$.
Hence  the connecting map for
the exact sequence $0 \to {\bf Z} \to {\bf Q} \to {\bf Q}/{\bf Z} \to 0$ yields an isomorphism $H^1(L,{\bf Q}/{\bf Z}) \simeq H^2(L,{\bf Z})$. Thus  $\sha^1_\omega(k,{\bf Q}/{\bf Z}) \simeq \sha^2_\omega(k, {\bf Z})$.
Since $\sha^1_\omega(k,{\bf Q}/{\bf Z})$  classifies   cyclic extensions of $k$ (together with a generator of the Galois group) which 
are locally almost everywhere split, by Chebotarev density theorem,  $\sha^1_\omega(k,{\bf Q}/{\bf Z}) = 0$.
In particular  $\sha^2(k,{\bf Z}) = 0$.

\begin{lemma}
\label{image} Let $M$ be a $\Gamma_k$-module and   $G_M$  the kernel of the map $\Gamma_k \to {\rm Aut}(M)$.
Let   $G = \Gamma_k/G_M $. The  image of the homomorphism  $H^2(G, M) \to H^2(k, M) $ contains $\sha^2_\omega(k, M)$.
\end{lemma}

\noindent
{\bf Proof.} Let $x \in \sha^2_\omega(k, M)$. Let $k_M = (k_s)^{G_M}$. Since $M$ becomes isomorphic to ${\bf Z}^n$ over 
$k_M$,   by Lemma \ref{Z} $x$ restricts to zero in $H^2(k_M, M)$.  Hence from the exact sequence $(*)$, $x$ belongs to the image of 
$H^2(G, M) \to H^2(k, M)$. 

\begin{lemma}
\label{Sansuc cycl} Let $M$ be a $\Gamma_k$-module,  let $G_M$ be the kernel of the map $\Gamma_k \to {\rm Aut}(M)$ and let
that $\Gamma_k/G_M = G$. Then $\sha_{\omega}^2(k,M) =  \sha^2_{\rm cycl}(k,M).$
\end{lemma} 

\noindent
{\bf Proof.} If $v \in \Omega_k$ is unramified in $k_M/k$, then $G_v$ is cyclic. Hence  
$\sha^2_{\rm cycl}(k,M) \subset \sha_{\omega}^2(k,M)$. 
 We show that $\sha_{\omega}^2(k,M) \subset  \sha^2_{\rm cycl}(k,M)$.  Let $x \in  \sha^2_\omega(k,M)$. 
 By Lemma \ref{image}, there is $y \in H^2(G,M)$ mapping to $x \in H^2(k,M)$.
 Let $v \in \Omega_k$ be such that $G_v$ is cyclic. 
By Chebotarev's density
theorem, there exist infinitely many $w \in \Omega_k$ such that $G_{w}$  is conjugate to  $G_{v}$.
Pick $w \in \Omega_k$ such that $x$ maps to zero in  $H^2(k_w, M)$ and there is a  $g \in G$  such that 
$gG_wg^{-1}=G_v$. The map $\psi_g:  M \to M$ given by $m \to gm$  is ${\rm Int}(g)$ semilinear and induces an isomorphism
$H^2(\psi_g): H^2(G, M) \to H^2(G, M)$ which is the identity \cite{HS} Proposition 16.2. Further, $H^2(\psi_g)$ restricts to  an isomorphism  $H^2(G_w, M) \to H^2(G_v, M)$.  Thus  the restriction of $y$ in $H^2(G_w, M)$ being  zero, its  restriction to $H^2(G_v, M)$ is zero and therefore  its image in $H^2(k_v, M)$ is zero. 
But this coincides with the image of $x$ in $H^2(k_v, M)$.
Thus $x$ maps to zero in $H^2(k_v, M)$. This is true for every  $v$ with $G_v$ is cyclic so that $x \in \sha^2_{\rm cycl}(k, M)$.

\bigskip
Let $\iota : T(k) \to \underset{v \in \Omega_k} \prod T(k_v)$ be the diagonal embedding, and let $A(T)$ be the quotient
of $\underset{v \in \Omega_k} \prod T(k_v)$  by the closure of the image of $\iota$; the group $A(T)$ is the obstruction to weak
approximation on $T$. Set $\sha(T) = \sha^1(k,T)$; this is the obstruction to the Hasse principle for torsors under $T$.

\bigskip
The following is a reformulation of a result of Voskresenskii :

\begin{prop}\label{Vosk} Let $G$ be a finite group, let $T$ be a $k$-torus, and assume that the character group of $T$ is a $G$-lattice via a surjection $\Gamma_k \to G$. We have
an exact sequence

$$0 \to A(T) \to \sha^2_{\rm cycl}(G,\hat T)^{\ast} \to \sha(T) \to 0.$$

\end{prop}

\medskip
\noindent
{\bf Proof.} Let $T^c$ be a smooth compactification of $T$;
by  \cite{San}, Theorem 9.5. (M) we have the exact sequence 
$$0 \to A(T) \to {\rm Br}_a(T^c)^{\ast} \to \sha(T) \to 0,$$  

Note that since $T_{k_s}$ is split and hence rational, we have ${\rm Br}(T^c_{k_s}) =0$
%${\rm Br}_a(T^c) = {\rm Br}(T^c)/{\rm Br}(\ell)$, since $T_{k_s}$ is rational, 
(see \cite {CTSk}, Corollary 5.2.6).
%hence ${\rm Br}(T^c_{k_s}) = 0$. By 
By Proposition \ref{Colliot Br cycl}
% \cite {CTS 87}, Theorem 9.5 (ii) 
we have ${\rm Br}(T^c)/{\rm Br}(k) \simeq \sha^2_{\rm cycl}(G,\hat T)$, 
hence we
 get the exact sequence 
$$0 \to A(T) \to \sha^2_{\rm cycl}(G,\hat T)^{\ast} \to \sha(T) \to 0.$$

\section{The group $\sha_{cycl}^2(G,M)$}\label{sha}

\medskip
Let $G$ be a finite group.

\begin{prop}\label{lemma} Let 
$\ell'/\ell$ be a finite Galois extension of global fields with Galois group $G$ which is unramified at all the finite places.
Let $M$ be a $G$-lattice regarded as a $\Gamma_{\ell}$-module via the surjection $\Gamma_{\ell} \to G $ . Then we have

$$\sha^2_{\rm cycl}(G,M)  \simeq  \sha^2(\ell,M).$$

\end{prop}

This proposition is an immediate consequence of Proposition \ref{omega} below :

\begin{prop}\label{omega} Let 
$\ell'/\ell$ be a finite Galois extension of global fields with Galois group $G$. 
Assume that all the decomposition groups of $\ell'/\ell$ are cyclic. Let $M$ be a $G$-lattice, regarded as $\Gamma_{\ell}$-lattice through
the surjection $\Gamma_{\ell} \to G$. Then we have

$$\sha^2_{\rm cycl}(G,M)  \simeq  \sha^2(\ell,M).$$

\end{prop}

Proposition \ref{omega} follows from Proposition \ref{cycl cycl} below. We use the notation of the previous section.

\begin{prop}\label {cycl cycl} Let $\ell$ be a global field, let $M$ be a $\Gamma_{\ell}$-module, and assume that $\Gamma_{\ell}/G_M \simeq G$. Then 
we have 

$$\sha^2_{\rm cycl}(G,M)  \simeq  \sha^2_{\rm cycl}(\ell,M).$$

\end{prop}

\medskip
\noindent
{\bf Proof of Proposition \ref{cycl cycl}.} Set $\ell' = \ell_M$; note that $\ell'/\ell$ is a Galois extension with group $G$.  The homomorphism $f : H^2(G,M) \to H^2(\ell,M)$ induced by the surjection $\Gamma_{\ell} \to G$ is
injective by Lemma \ref{1.3}. By Lemmas \ref{image} and  \ref{Sansuc cycl}, 
 the image of $f$ contains $\sha^2{\rm cycl}(\ell, M)$.
We next prove that the restriction of $f$ to $\sha^2_{\rm cycl}(G,M)$ maps it into  $\sha^2_{\rm cycl}(\ell,M)$. Let $x \in \sha^2_{\rm cycl}(G,M)$. Let $v \in \Omega_{\ell}$
such that 
%and let $G_v$ be the decomposition group of $\ell'/\ell$ at $v$. By hypothesis, 
 the decomposition group $G_v$ of $\ell'/\ell$ at $v$ is a cyclic subgroup of $G$. Then the restriction of 
$x \in H^2(G,M)$ to $H^2(G_v,M)$ is zero. The composite $H^2(G,M) \to H^2(\ell,M) \to H^2(\ell_v,M)$ factors as 
$H^2(G,M) \to H^2(G_v,M) \to H^2(\ell_v,M)$. Hence $f(x)$ maps to zero in $H^2(\ell_v,M)$. Thus  $f(x) \in \sha^2_{\rm cycl}(\ell,M)$.

\medskip
Clearly $\sha^2_{\rm cycl}(G,M)  \to  \sha^2_{\rm cycl}(\ell,M)$ is injective. We prove that this map is surjective. Let $y \in \sha_{\rm cycl}^2(\ell,M)$ and
let $x \in H^2(G,M)$ be such that $f(x) = y$. Let $g \in G$. 
By Chebotarev's density theorem, there is a finite place $v \in \Omega_{\ell}$  such that $G_v = \langle g \rangle$.
We claim that the restriction of $x$ to $H^2(\langle g \rangle, M) = H^2(G_v, M)$ maps to zero in 
$H^2(\ell_v, M)$. 
In fact this image is the same as the restriction of $y$ to $H^2(\ell_v, M)$.
Since $ y \in \sha^2_{\rm cycl}(\ell, M)$ and $G_v$ is cyclic,  the image of $y$ in $H^2(\ell_v, M)$
is zero. 
By lemma \ref{1.3}, the map $H^2(G_v,M) \to H^2(\ell_v,M)$ is injective. It follows that the restriction of $x$ to $H^2(G_v,M)$ is zero and
hence $x$ belongs to $\sha^2_{\rm cycl}(G,M)$.
This completes the proof of the proposition.

\medskip
\noindent
{\bf Proof of Proposition \ref {omega}.} Since all the decomposition groups of $\ell'/\ell$ are cyclic, we have $ \sha^2(\ell,M) =  \sha^2_{\rm cycl}(\ell,M).$
By Proposition \ref {cycl cycl}, we have $\sha^2_{cycl}(G,M)  \simeq  \sha^2_{\rm cycl}(\ell,M)$, hence 
$\sha^2_{cycl}(G,M)  \simeq  \sha^2(\ell,M)$, as claimed.

\medskip
\noindent
{\bf Proof of Proposition \ref {lemma}.} Since $\ell'/\ell$ is unramified at all the finite places, all the decomposition groups are cyclic; we conclude by
applying Proposition \ref {omega}.

\begin{coro}\label{Sansuc} Let $\ell$ be a global field, let $M$ be a $\Gamma_{\ell}$-module, and assume that $\Gamma_{\ell}/G_M \simeq G$. Then 
we have 
Theorem \ref{p+2} (b)
$$\sha^2_{\rm cycl}(G,M)  \simeq  \sha^2_{\omega}(\ell,M).$$

\end{coro}

\medskip
\noindent
{\bf Proof.} This follows from Proposition \ref{cycl cycl} and Lemma \ref {Sansuc cycl}.

\begin{coro}\label{cycl F} Let $M$ be a $G$-lattice. Let $\ell'/\ell$ be as in Proposition \ref {omega}, and let $T$ be an $\ell$-torus with character group $M$. We have

$$\sha^2_{\rm cycl}(G,M)  \simeq \sha^2(\ell,M) = \sha^2(\ell,\hat T) \simeq \sha^1(\ell,T)^{\ast}.$$

\end{coro}

\medskip
\noindent{\bf Proof.} This follows from Proposition \ref{omega}, and from Poitou-Tate duality.

\medskip
In the following sections, we also need a result of Fr\"ohlich :

%\medskip
%The following results are well-known

\begin{prop}\label{Frohlich} {\rm (Fr\"ohlich)} There exists a Galois extension of number fields with Galois group $G$ that is unramified at all the finite places.

\end{prop}

\noindent
{\bf Proof.} This is the main result of \cite{F}.

%\begin{prop}\label{equivalence of categories} Let $F$ be a field, and let $M$ be a $G$-lattice. There exists an $F$-torus $T$ such that the character
%lattice of $T$ is isomorphic to the $G$-lattice $M$, regarded as $\Gamma_F$-lattice through the surjection $\Gamma_F \to G$. 

%\end{prop}

%\noindent 
%{\bf Proof.} See \cite {Bo}, Chapter III, \S 8.

\medskip
In the next sections we use Corollary \ref{cycl F} together with \ref{Frohlich} to realize the purely algebraic group $\sha^2_{\rm cycl}(G,M)$
as the Tate-Shafarevich group of a torus over a number field. This makes it possible to apply arithmetic results to obtain algebraic ones.

\section{Vanishing results}

The aim of this section is to apply the results of \S \ref{sha} and, under an additional hypothesis (condition (C) below) prove some vanishing theorems.
Let $G$ be a finite group.

\begin{defn}\label{C}
Let $M$ be a $G$-lattice. We say that $M$ {\it satisfies condition {\rm (C)}} if there exists
a Galois extension $\ell'/\ell$ of number fields with Galois group $G$ such that all the decomposition groups of $\ell'/\ell$ are cyclic,
and such that the  $\ell$-torus $S$ associated to the Galois lattice $M$ with the Galois group $\Gamma_{\ell}$ acting through the quotient group G has the property $\sha^1(\ell, S) = 0$.

\end{defn}

\begin{coro}\label{0} Let $M$ be a $G$-lattice satisfying condition  {\rm (C)}. 
% Let $\ell'/\ell$ be as in Proposition \ref {omega}, and let $S$ be an $\ell$-torus with
%character $G$-lattice $M$. Assume that $\sha^1(\ell,S) = 0$. 
Then $$\sha^2_{\rm cycl}(G,M) = 0.$$

\end{coro}

\noindent
{\bf Proof.} 
By Poitou-Tate duality, we have $\sha^1(\ell,S) \simeq \sha^2(\ell,\hat S)^{\ast} = \sha^2(\ell,M)^{\ast}$.  Corollary
 \ref{cycl F} implies that the groups $\sha^2(\ell,\hat S) = \sha^2(\ell,M)$ and $\sha^2_{cycl}(G,M)$ are
isomorphic. Hence $\sha^1(\ell,S)$ is dual to the group  $\sha^2_{cycl}(G,M)$. 
Since   
$\sha^1(\ell,S) = 0$, we have $\sha^2_{cycl}(G,\hat S) =  \sha^2_{cycl}(G,M)  =0$.

\begin{coro}\label{hp wa} Let $k$ be a global field with a surjection $\Gamma_k \to G$, and let $T$ be a $k$-torus; assume that the character group of $T$ is a $G$-lattice satisfying  condition  {\rm (C)}.
%satisfies the hypotheses of 
%Corollary \ref{0}. 
Then 
$$\sha^2_{\omega}(k,\hat T)= 0,$$
and  Hasse principle and weak approximation hold for torsors under $T$. 

\end{coro}

\noindent
{\bf Proof.} %Let $\iota : T(k) \to \underset{v \in \Omega_k} \prod T(k_v)$ be the diagonal embedding, and let $A(T)$ be the quotient
%of $\underset{v \in \Omega_k} \prod T(k_v)$  by the closure of the image of $\iota$; the group $A(T)$ is the obstruction to weak
%approximation on $T$. Set $\sha(T) = \sha^1(k,T)$; this is the obstruction to the Hasse principle for torsors under $T$. 
%By  \cite{San}, Theorem 9.5. (M) have the exact sequence 
%$$0 \to A(T) \to {\rm Br}_a(T)^{\ast} \to \sha(T) \to 0,$$  
%Note that ${\rm Br}_a(T) = {\rm Br}(T)/{\rm Br}(\ell)$, since $T_{k_s}$ is rational, hence ${\rm Br}(T_{k_s}) = 0$. By 
 %\cite {CTS 87}, Theorem 9.5. (ii) we have ${\rm Br}(T)/{\rm Br}(\ell) \simeq \sha^2_{\rm cycl}(G,\hat T)$, hence we
 %get the exact sequence 
 By Proposition \ref{Vosk}, we have the exact sequence
$$0 \to A(T) \to \sha^2_{\rm cycl}(G,\hat T)^{\ast} \to \sha(T) \to 0,$$  We have
$\sha^2_{\rm cycl}(G,\hat T)^{\ast} = 0$ by Corollary \ref{0}, hence $A(T) = \sha(T) = 0$. By Corollary \ref{Sansuc}, we have  
$\sha_{\omega}(k,\hat T)= 0$; weak approximation holds
for $T$, and Hasse principle holds for torsors under $T$.

\medskip This implies Colliot-Th\'el\`ene's observation cited in the introduction :

\begin{coro}\label{all fields}   Let $G$ be a finite group, and let $M$ be a $G$-lattice. If for all number fields $k$ and every $k$-torus $T$ of character group isomorphic to
the Galois module $M$ via a surjection of $ \Gamma_k \to G$, one can show that $\sha^1(k,T) = 0$, then $\sha_{\rm cycl}^2(G,M) = 0$, and  every such $k$-torus $T$ satisfies weak approximation.
\end{coro}

%\begin{remark}\label{decomposition groups} Note that the hypothesis on $\ell'/\ell$ of Proposition \ref {lemma} and Corollaries \ref {cycl F} and \ref {0} can be weakened; it is enough to assume that all the decomposition groups
%of $\ell'/\ell$ are cyclic. The proofs are the same.

%\end{remark}

\section{Unramified Brauer groups}\label{un}

Let $k$ be a field,  $T$  a $k$-torus  and  $X$  a torsor under $T$. Let $T^c$ be a smooth equivariant compactification of $T$; such a compactification exists in any characteristic (see \cite{CTHSk 05}). Then, the contracted product  $X^c = X \times^TT^c$  is a smooth compactification of $X$.
Note that $\overline T^c \simeq \overline X^c$ is birational to the
projective space; further $T^c(k) \neq \emptyset$. We therefore have exact sequences,  by Theorem  \ref{brnr},

$${\rm Br}(k) \to {\rm Br}(X^c) \to H^1(k,{\rm Pic}(\overline X^c))$$ 
and 
$$0 \to {\rm Br}(k) \to {\rm Br}(T^c) \to H^1(k,{\rm Pic}(\overline T^c)) \to 0.$$ 

\medskip
By \cite[Lemma 2.1]{CTHSk 03}, $H^1(k,{\rm Pic}(\overline X^c)) \simeq H^1(k,{\rm Pic}(\overline T^c))$, we have an injection

$${\rm Br}(X^c)/{\rm Im}({\rm Br}(k)) \to {\rm Br}(T^c)/{\rm Br}(k).$$

%The {\it unramified Brauer group} of $X$ is by definition the group ${\rm Br}(X^c)$.
%the quotient ${\rm Br}(X^c)/{\rm Br}(k)$. 
\medskip
\noindent
The aim of this section and the following ones is to use the results of the previous sections to obtain information on the quotients ${\rm Br}(T^c)/{\rm Br}(k) $ and
${\rm Br}(X^c)/{\rm Im}({\rm Br}(k) )$. We start with some vanishing results, 

\bigskip
%\noindent
%Let $G$ be a finite group. Recall that a $G$-lattice $M$  satisfies condition {\rm (C)}  if there exists
%a Galois extension $\ell'/\ell$ of number fields with Galois group $G$ such that all the decomposition groups of $\ell'/\ell$ %are cyclic,
%and such that the the $\ell$-torus $S$ associated to the Galois lattice $M$ with the Galois group $\Gamma_{\ell}$ acting %through the quotient group G has the property $\sha^1(\ell, S) = 0$ (cf. definition \ref{C}). 

\begin{prop}\label{general} Let $G$ be a finite group.  Assume that the character  group of $T$  is a $G$-lattice satisfying  condition  {\rm (C)}.
%Corollary \ref{0}. 
Then ${\rm Br}(T^c)/{\rm Br}(k) = 0.$ If $X$ is a torsor over $T$, we have  
${\rm Br}(X^c)/{\rm Im}({\rm Br}(k)) = 0$.
\end{prop}

\medskip
\noindent
{\bf Proof.} By Corollary \ref{0}, we have $\sha^2_{\rm cycl}(G,\hat T) = 0$. 
On the other hand,  by Theorem  \ref{Colliot Br cycl},
  ${\rm Br}(T^c)/{\rm Br}(k) \simeq
\sha^2_{\rm cycl}(G,\hat T)$ ), therefore ${\rm Br}(T^c)/{\rm Br}(k) = 0$. As we saw above,
${\rm Br}(X^c)/{\rm Im}({\rm Br}(k))$ injects into ${\rm Br}(T^c)/{\rm Br}(k)$, hence this implies that ${\rm Br}(X^c)/{\rm Im}({\rm Br}(k)) = 0$.

\section{Norm equations}\label{norm}
 
\medskip
Let $k$ be a field, and let $L$ be an \'etale $k$-algebra of finite rank (in other words, a product of a finite number of
separable extensions of $k$). Let $T_{L/{k}}  = R^{(1)}_{L/k}({\bf G}_m)$ be the $k$-torus
defined by $$1 \to T_{L/{k}} \to R_{L/k}({\bf G}_m) {\buildrel {\rm N}_{L/{k}} \over \longrightarrow}  {\bf G}_m \to 1.$$

\medskip Let $a \in k^{\times}$. Let $X$ be the affine $k$-variety associated to the {\it norm equation}

$${\rm N}_{L/k}(t) = a.$$

\medskip
The variety $X$ is a torsor under $T_{L/k}$;
let  $X^c$ be a smooth compactification of $X$. 

\medskip
In the later sections, we need the following result

\begin{theo}\label{T = X} Suppose $L=K \times E$ with $K/k$  a cyclic extension  and  $E$  an \'etale $k$-algebra of finite rank. Then

$${\rm Br}(T^c)  / {\rm Br}(k) \simeq 
{\rm Br}(X^c)/{\rm Im}({\rm Br}(k)).$$

\end{theo}

\noindent
{\bf Proof.}
%Let $X$ be defined by the equation $N_{L/k}(t) = c$, $c \in k^*$. 
Since $K\otimes_kL \simeq K \times L'$ for some \'etale algebra $L'$ over $K$, the variety
$X_K$ is isomorphic to $R_{L'/K}(G_m)$  and hence is $K$-rational. Since $K/k$ is cyclic, by \cite[Proposition 1.1]{CT}, the map 
$H^1(k, {\rm Pic}\overline{X}^c) \to H^3(k, {\bf G}_m)$ in the sequence $(**)$ is zero and one has an isomorphism 
(cf. Theorem \ref{brnr}) ${\rm Br}(X^c) /{\rm Im} ( {\rm Br}(k)) \simeq H^1(k, {\rm Pic}\overline{X}^c)$. We also have an isomorphism (cf. Theorem \ref{brnr})
$${\rm Br}(T^c)/{\rm Br}(k) \simeq H^1(k, {\rm Pic}\overline{T}^c).$$
By \cite[Lemma 2.1]{CTHSk 03}, we have an isomorphism 
$$H^1(k, {\rm Pic}\overline{X}^c) \simeq H^1(k, {\rm Pic}\overline{T}^c).$$
This leads to an isomorphism 
$${\rm Br}(T^c)  / {\rm Br}(k) \simeq 
{\rm Br}(X^c)/{\rm Im}({\rm Br}(k)).$$

\section{Norm equations - first examples}\label{norm examples}

The aim of this section and the next ones is to give some examples of \'etale algebras for which we apply the results of the previous sections,
obtaining information about the unramified Brauer group,  Hasse principle and weak approximation 
(in the case where $k$ is a global field) for the variety $X$. We keep the notation of the previous section.

\medskip
The first examples concern \'etale algebras that are products of two fields, finite extensions of the ground field $k$.

\medskip
{\bf Products  of two fields}

\medskip We start by introducing some notation that will be used in the two examples of this section. Let $L = K_1 \times K_2$, where $K_1$ and $K_2$ are finite extensions of $k$. Let $k'/k$ be a Galois extension of minimal degree splitting $T_{L/{k}}$, and let $G = {\rm Gal}(k'/k)$; let
$M = \hat T_{L/k}$ be the character $G$-lattice of $T_{L/k}$. For $i = 1,2$, let $H_i$ be the subgroup of $G$ such that $K_i = (k')^{H_i}$.We have the exact sequence
of $G$-modules

$$0 \to {\bf Z} \to  {\bf Z}[G/H_1]  \oplus {\bf Z}[G/H_2] \to M \to 0.$$

\medskip
Let  $\ell'/\ell$ be an unramified extension of number fields with Galois group $G$.

\medskip
{\bf Some consequences of H\"urlimann's theorem}

\medskip The first example is based on a result of H\"urlimann \cite{H}.  With the notation above, we assume that $K_1/k$ is a cyclic extension. 

\begin{theo} \label {cyclic Hurlimann} We have $\sha^2_{\rm cycl}(G,M) = 0$. 

\end {theo}

\noindent
{\bf Proof.} Set $\ell_i = (\ell')^{H_i}$ 
for $i = 1,2$. Let $S$ be the norm torus corresponding to the \'etale $\ell$-algebra $\ell_1 \times \ell_2$. H\"urlimann's result \cite {H}, Proposition 3.3  implies that
$\sha^1(\ell,S) = 0$ (in \cite {H}, the extension $K_2/k$ is supposed to be Galois, but this is not necessary; see \cite {BLP}, Proposition 4.1  for a different proof of
the general case). We have $\hat S \simeq M$ by construction, hence by Proposition \ref{lemma} we have $\sha^2_{\rm cycl}(G,M) = 0$. 

\begin{remark} Theorem \ref {cyclic Hurlimann} was proved by Sansuc (unpublished) by algebraic methods. His proof is rather involved; the
proof presented here, passing from arithmetic to algebra,  is simpler, since the proof of the arithmetic result in  \cite {BLP}, Proposition 4.1 is quite short.

\end{remark}

\begin{theo} \label {Hurlimann}   We have ${\rm Br}(T_{L/k}^c)  / {\rm Br}(k) = 0$, and ${\rm Br}(X^c)  / {\rm Im}({\rm Br}(k)) = 0$.

\end{theo} 

\noindent
{\bf Proof.} By Theorem \ref{cyclic Hurlimann} and Proposition \ref{general}, we have ${\rm Br}(T_{L/k}^c)  / {\rm Br}(k) = 0$ 
and  ${\rm Br}(X^c)  / {\rm Im}({\rm Br}(k)) = 0$.

\bigskip
{\bf Linearly disjoint Galois extensions}

\medskip
This example is based on a result of Pollio and Rapinchuk, \cite{PR}. Let
$L = K_1 \times K_2$, where $K_1$ and $K_2$ are finite extensions of $k$ such that the Galois closures of $K_1$ and $K_2$ are linearly disjoint.

\begin{theo} \label {Polliot-Rapinchuk}   We have ${\rm Br}(T_{L/k}^c)  / {\rm Br}(k) = 0$, and ${\rm Br}(X^c)  / {\rm Im}({\rm Br}(k)) = 0$.

\end{theo}

\noindent 
{\bf Proof.} %Let $T = T_{L/k}$ as above. 
 Let $K_1'$ and $K_2'$ be the Galois closures of $K_1$, respectively $K_2$,
and let $H_1'$, $H_2'$ be the subgroups of $G$ such that $K'_i = (k')^{H'_i}$, for $i = 1,2$. 
By hypothesis, the extensions $K_1'$ and $K_2'$ are linearly disjoint, hence  $G = H'_1.H'_2$.

\medskip
Recall that   $\ell'/\ell$ is an unramified extension of number fields with Galois group $G$. Set $\ell_i = (\ell')^{H_i}$ and $\ell'_i = (\ell')^{H'_i}$
for $i = 1,2$. Since $H_i'$ is a normal subgroup of $G$ for $i = 1,2$, the extensions 
$\ell_i'/\ell$ are Galois, and 
$|H_i'| = [\ell': \ell_i']$. Since
$G = H'_1.H_2'$, the fields $\ell_i'$ are linearly disjoint.  

\medskip
Let $S$ be the norm torus corresponding to the \'etale $\ell$-algebra $\ell_1 \times \ell_2$; the main theorem of \cite{PR} implies that
$\sha^1(\ell,S) = 0$. We have $\hat S \simeq M$ by construction, hence by Corollary \ref{0}  we have $\sha^2_{\rm cycl}(G,M) = 0$. 
By Proposition \ref{general}, we have ${\rm Br}(T_{L/k}^c)  / {\rm Br}(k) = 0$ and ${\rm Br}(X^c)  / {\rm Im}({\rm Br}(k)) = 0$.

\section{Products of cyclic extensions of prime power degree - statement of results and notation}\label{statement}

\medskip The proofs of the results of this section will be given in \S \ref {proof}.
Let $p$ be a prime number.
If $K/k$ is a cyclic extension of degree a power of $p$, we denote by $(K)_{\rm prim}$ the unique 
subfield of $K$ of degree $p$ over $k$; if $E = \underset {i \in I} \prod K_i$, where $K_i/k$ is a cyclic extension of
degree a power of $p$ for all $i \in I$, set $E_{\rm prim} = \underset {i \in I} \prod (K_i)_{\rm prim}$. 

\medskip
Let $L$ be a product of $n$ cyclic extensions of degrees powers of $p$. With the notation of the
previous section, set $T = T_{L/k}$.
%and $T_{\rm prim} = T_{L_{\rm prim}/k}$. 
Let $a \in k^{\times}$, and $X$ be the affine $k$-variety associated to the norm equation
${\rm N}_{L/k}(t) = a$; 
let  $X^c$ be a smooth compactification of $X$.  

%Similarly, we denote by $X_{\rm prim}$ the affine $k$-variety associated to ${\rm N}_{L_{\rm prim}/k}(t) = a$,
%and by $X^c_{\rm prim}$ a smooth compactification of $X_{\rm prim}$.  
\medskip

Let $T_{\rm prim} = T_{L_{\rm prim}/k}$. We denote by $X_{\rm prim}$ the affine $k$-variety associated to ${\rm N}_{L_{\rm prim}/k}(t) = a$,
and by $X^c_{\rm prim}$ a smooth compactification of $X_{\rm prim}$. 

\medskip

\begin{theo}\label{torus p} We have $${\rm Br}(T^c)  / {\rm Br}(k) = 0 
\iff  {\rm Br}(T_{\rm prim}^c)  / {\rm Br}(k) = 0.$$

\end{theo}

\medskip
Let $k'/k$ be a Galois extension of minimal degree splitting $T$, and let $G = {\rm Gal}(k'/k)$; similarly, let
$k_{\rm prim}'/k$ be a Galois extension of minimal degree splitting $T_{\rm prim}$, and let $G_{\rm prim} = {\rm Gal}(k'_{\rm prim}/k)$. 

\medskip

\begin{theo}\label{cyclic p} We have $$\sha^2_{\rm cycl}(G,\hat T) = 0 \iff \sha^2_{\rm cycl}(G_{\rm prim},\hat T_{\rm prim}) = 0.$$

\end{theo}

\bigskip

{\bf Products of at least $p+2$ pairwise disjoint cyclic extensions}

\bigskip
With the notation above, we now consider the case where $n \ge p+2$.

\begin{theo} \label {p+2}  Assume that $L$ is a product of at least
$p+2$ pairwise disjoint cyclic extensions of degrees powers of $p$. 
 Then we have

\medskip
{\rm (a)}  $${\rm Br}(T^c)  / {\rm Br}(k) = 0.$$

\medskip
{\rm (b)}  $${\rm Br}(X^c)  /{\rm Im}( {\rm Br}(k) )= 0.$$

\medskip
{\rm (c)} Suppose that $k$ is a global field. Then $\sha^2_{\omega}(k,\hat T_{L/k})= 0$,  and Hasse principle and weak approximation hold for $X$.

\end{theo}

\smallskip

\begin{theo} \label {p+2 cycl} Assume that $n \ge p + 2$. Then we have
$$\sha^2_{\rm cycl}(G,\hat T) = 0.$$

\end{theo}

\bigskip

{\bf At least one cyclic factor of degree $p$}

\bigskip

In the next results, we assume that $L$ has at least one factor of degree $p$.

%\medskip
%Let $L$ be a product of distinct cyclic extensions of degree a power of $p$; assume moreover that {\it at
%least one of these is of degree $p$}. With the notation of the previous section, set $T = T_{L/k}$ and $T_{\rm prim} = T_{L_{\rm prim}/k}$.

\medskip

\begin{theo}\label{torus at least one prime} Assume that $L$ is a product of  $n$ pairwise
disjoint cyclic extensions of degrees powers of $p$, and that at least one of these is of degree $p$.  Then we have $${\rm Br}(T^c)  / {\rm Br}(k)  
\simeq {\rm Br}(T_{\rm prim}^c)  / {\rm Br}(k).$$

\end{theo}

\medskip

%Let $k_{\rm prim}'/k$ be a Galois extension of minimal degree splitting $T_{\rm prim}$, and let $G_{\rm prim} = {\rm Gal}(k'_{\rm prim}/k)$. 

\begin{theo}\label{cyclic at least one prime} Assume that $L$ is a product of $n$ pairwise
disjoint cyclic extensions of degrees powers of $p$, and that at least one of these is of degree $p$. Then we have $$\sha^2_{\rm cycl}(G,\hat T)  \simeq \sha^2_{\rm cycl}(G_{\rm prim},\hat T_{\rm prim}).$$

\end{theo}

\bigskip
{\bf Products of cyclic extensions of degree $p$}

\bigskip

Finally, we determine ${\rm Br}(T_{\rm prim}^c)  / {\rm Br}(k)$. Let us denote by $C_p$ the cyclic group of order $p$.

\begin{theo}\label{prime}

{\rm (a)} If $G_{\rm prim} \not \simeq C_p \times C_p$, then $${\rm Br}(T_{\rm prim}^c)  / {\rm Br}(k) =  {\rm Br}(X_{\rm prim}^c)  / {\rm Im}({\rm Br}(k) ) = 0.$$

\medskip
{\rm (b)} If $G_{\rm prim}  \simeq C_p \times C_p$, then 
$${\rm Br}(T_{\rm prim}^c)/{\rm Br}(k)  \simeq {\rm Br}(X_{\rm prim}^c)  / {\rm Im}({\rm Br}(k) ) \simeq ({\bf Z} / p {\bf Z})^{n-2}.$$

\end{theo}

\begin{theo}\label{prime cyclic sha}

{\rm (a)} If $G_{\rm prim} \not \simeq C_p \times C_p$, then $\sha^2_{\rm cycl}(G_{\rm prim},\hat T_{\rm prim}) = 0$.

\medskip
{\rm (b)} If $G_{\rm prim}  \simeq C_p \times C_p$, then $\sha^2_{\rm cycl}(G_{\rm prim},\hat T_{\rm prim}) \simeq ({\bf Z} / p {\bf Z})^{n-2}.$

\end{theo}

\medskip In the case where $k$ is a global field, we obtain the following corollaries

\begin{coro}\label{weak approx} Suppose that $k$ is a global field. Assume that $L$ is a product of  $n$ pairwise
disjoint cyclic extensions of degrees powers of $p$, and that at least one of these is of degree $p$.  

\medskip

{\rm (a)} If $G_{\rm prim} \not \simeq C_p \times C_p$, then Hasse principle and weak approximation hold for $T $.

\medskip
{\rm (b)} If $G_{\rm prim}  \simeq C_p \times C_p$, then either Hasse principle holds for torsors over $T$ (and weak approximation for
$T$ fails), 
or 
weak approximation holds for $T$ (and Hasse principle for torsors over $T$ fails). 

\end{coro}

\noindent
{\bf Proof.} Part (a) follows from Theorems \ref{prime cyclic sha} (a), 
\ref {cyclic at least one prime} and  \ref {Vosk}. To prove part (b), we apply Theorems \ref{prime cyclic sha} (b), \ref {cyclic at least one prime},   \ref {Vosk},
as well as \cite {BLP}, Theorem 8.3 and Corollary 5.17.

\begin{coro}\label{corocoroproduct} Assume that $k$ is a global field, and that $L$ is a product of $n$ distinct cyclic extensions of degree $p$. 
If $G \not \simeq C_p \times C_p$, then $$\sha_{\omega}(k,\hat T_{L/k}) = 0.$$ If $G \simeq C_p \times C_p$, then $$\sha_{\omega}(k,\hat T_{L/k}) \simeq 
({\bf Z} / p {\bf Z})^{n-2}.$$

\end{coro}

\noindent
{\bf Proof.} This follows from Theorem \ref{prime cyclic sha} and Lemma \ref{Sansuc}.

\medskip
\noindent
\begin{remark}\label{macedo}
Corollary \ref{corocoroproduct} was also proved by Macedo, see \cite {Ma}, Theorem 4.9 and Corollary 4.10, with different methods
(namely, using a generalization of the approach of Drakokhrust and Platonov).

\end{remark}

\medskip

\section{Products of cyclic extensions of prime power degree - global fields}\label{product global}

We recall some results from \cite{BLP}; these will be used in the next section to prove the results of \S \ref{statement}. Let $k$ be a global field. We start by recalling some notation from \cite {BLP}. If $L$ is an
\'etale algebra of finite rank over $k$ having at least one factor that is a cyclic extension of $k$,  the paper \cite {BLP} introduces a finite abelian 
group $\sha(L)$ (see \cite {BLP},  \S 5) and proves (see \cite {BLP}, Corollary 5.17) that $\sha(L)^{\ast} \simeq \sha^1(k,T_{L/k})$; equivalently, 
by Poitou-Tate duality, we have $\sha (L) \simeq \sha^2(k,\hat T_{L/k})$.

\medskip
Let $p$ be a prime number. 

\begin{prop}\label{injection} Let $L$ be a product of cyclic extensions of degrees powers of $p$. The group $\sha^1(k,T_{L_{\rm prim}/k})^{\ast}$ injects into
$\sha^1(k,T_{L/k})^{\ast}$.

\end{prop}

\noindent
{\bf Proof.} This follows from \cite {BLP}, Lemma 8.7.

\begin{theo}\label{global}  Let $L$ be a product of cyclic extensions of degrees powers of $p$. Then we have $$\sha^1(k,T_{L/k}) = 0 \iff
 \sha^1(k,T_{L_{\rm prim}/k}) = 0.$$

\end{theo}

\noindent
{\bf Proof.} This is an immediate consequence of \cite {BLP}, Theorem 8.1. 

\begin{prop}\label {prime sha} Let $L$ be a product of $n$ distinct cyclic extensions of degrees powers of $p$, and assume that at least one of the  extensions is of degree $p$.
Then $\sha^1(k,T_{L/k})$ is a finite abelian group of type $(p,...,p)$ of order at most $p^{n-2}$.

\end{prop}

\noindent
{\bf Proof.} Let us write $L$ as a product $L = K \times K'$, where $K$ is a cyclic extension of $k$ of degree $p$, and $K'$ is a product of $n -1$ cyclic extensions of degrees
powers of $p$. 
In  \cite {BLP}, 5.1, we construct a finite abelian  group $\sha(L)=\sha(K,K')$ such that when $K$ is cyclic of
order $p$, the group $\sha(K,K')$ is of type $(p,...,p)$ of order at most $p^{n-2}$. It is shown in 
 \cite {BLP},  5.3 that the group $\sha(K,K')$
does not depend on the decomposition of $L$ as $K \times K'$, where $K$ is a cyclic extension of $k$, and that $\sha(K, K')^{\ast} \simeq \sha^1(k,T_{L/k})$
(see \cite {BLP}, Corollary 5.17). Hence $\sha^1(k,T_{L/k})$ is a finite abelian group of type $(p,...,p)$ of order at most $p^{n-2}$, as claimed.
\medskip

\begin{theo}\label{BLP} Let $L = K_1 \times \dots \times K_n$, where $K_i$ are distinct cyclic extensions of degree $p$ of $k$.

\medskip
 $(a)$ If $n \le 2$, or $n \ge p+2$, or $3 \le n \le  p+1$ and $K_1, \dots, K_n$  are not contained in some field extension 
of $k$
of degree $p^2$ having all local degrees $\le p$, then 
$$\sha^1(k,T_{L/{k}}) = 0.$$

\medskip $(b)$ Assume that $3 \le n \le p+1$ and that the fields $K_1, \dots, K_n$  are contained in some field extension of ${k}$
of degree $p^2$ having all local degrees $\le p$, then $$\sha^1(k,T_{L/{k}})  \simeq  ({\bf Z} / p {\bf Z})^{n-2}.$$

\end{theo}

\medskip
\noindent
{\bf Proof.} This follows from  \cite{BLP}, Theorem 8.3, Proposition 8.5 and Corollary 5.17. 

\medskip
\noindent
{\bf Remark.} More generally, one can treat the case where the $K_i$'s are field extensions of degree $p$, with at least one of them cyclic
(see \cite{BLP}, Proposition 8.5). 

\begin{theo}\label{global at least one prime} Let $L$ be a product of $n$ pairwise disjoint cyclic extensions of degrees powers of $p$, and assume that at least one of the  extensions is of degree $p$.
Then we have $$\sha^1(k,T_{L/k}) \simeq
 \sha^1(k,T_{L_{\rm prim}/k}).$$

\end{theo}

\noindent
{\bf Proof.} By Theorem \ref {BLP}, we have either $ \sha^1(k,T_{L_{\rm prim}/k}) = 0$ or $ \sha^1(k,T_{L_{\rm prim}/k})  \simeq  ({\bf Z} / p {\bf Z})^{n-2}.$
If $ \sha^1(k,T_{L_{\rm prim}/k}) = 0$, then by Theorem \ref {global}, we have $\sha^1(k,T_{L/k}) = 0$. Assume now that $ \sha^1(k,T_{L_{\rm prim}/k}) 
\not = 0$; then by Theorem \ref {BLP} we have $ \sha^1(k,T_{L_{\rm prim}/k})  \simeq  ({\bf Z} / p {\bf Z})^{n-2}$, therefore the order  of 
$ \sha^1(k,T_{L_{\rm prim}/k}) $ is equal to $p^{n-2}$. By Proposition \ref {injection}, this implies that the order of $\sha^1(k,T_{L/k})$ is at
least $p^{n-2}$. On the other hand, since at least one of the factors of $L$ is of order $p$, Proposition \ref {prime sha} implies that
$\sha^1(k,T_{L/k})$  is a finite abelian group of type $(p,\dots,p)$ of order
at most $p^{n-2}$. Hence the order of $\sha^1(k,T_{L/k})$ is equal to $p^{n-2}$, and this completes the proof of the Theorem.

\bigskip

\section{Products of cyclic extensions of prime power degree - proofs}\label{proof}

We keep the notation of \S \ref {statement}. In particular, $p$ is a prime number, $L$ is a product of $n$ cyclic extensions of degrees powers of $p$,
 $k'/k$ is a Galois extension of minimal degree splitting $T = T_{L/{k}}$, and  $G = {\rm Gal}(k'/k)$. Let $M = \hat T$ be the $G$-lattice of characters of $T$. Let us write $L = \underset {i \in I} \prod K_i$. Since $k'$ splits $T$, $k'$ also splits $R_{L/k}(G_m)$ and it follows that $k'$ contains all the factors $K_{i}$ of $L$.
 Let $H_i$ be the subgroup of $G$ such that $K_i = (k')^{H_i}$. We have the exact sequence
of $G$-modules

$$0 \to {\bf Z} \to \underset {i \in I} \oplus {\bf Z}[G/H_i] \to M \to 0.$$

Recall that $k_{\rm prim}'/k$ is a Galois extension of minimal degree splitting $T_{\rm prim}$, and that $G_{\rm prim} = {\rm Gal}(k'_{\rm prim}/k)$. 
Let $M_{\rm prim}$ be the $G_{\rm prim}$-lattice of characters of $T_{\rm prim}$.
Let $(H_i)_{\rm prim}$ be the subgroup of $G_{\rm prim}$ such that $(K_i)_{\rm prim} = (k_{\rm prim}')^{(H_i)_{\rm prim}}$. We have the exact sequence
of $G_{\rm prim}$-modules

$$0 \to {\bf Z} \to \underset {i \in I} \oplus {\bf Z}[G_{\rm prim}/(H_i)_{\rm prim}] \to M_{\rm prim} \to 0.$$

Set $\Gamma_i = (G/H_i)/(G_{\rm prim}/(H_i)_{\rm prim})$, and note that $(K_i)_{\rm prim} = K_i^{\Gamma_i}$.

\medskip
Let $\ell'/\ell$ be an extension of number fields   
with Galois group $G$ which is unramified at all the finite places.
For all $i \in I$, let $L_i$ be the fixed field of $H_i$ in $\ell'$, and set $E = \underset{i \in I} \prod L_i$. The character lattice 
of the torus $T_{E/\ell}$ is isomorphic to the $G$-lattice $M$. We have $(L_i)_{\rm prim} = L_i^{\Gamma_i}$.
%Set $\Gamma_i = (G/H_i)/(G_{\rm prim}/(H_i)_{\rm prim})$.
%Recall that for all $i \in I$, we denote by $(L_i)_{\rm prim}$ the unique 
%subfield of $L_i$ of degree $p$ over $k$, and we set $E_{\rm prim} = \underset {i \in I} \prod (L_i)_{\rm prim}$.
Note that $\sha^2_{\rm cycl}(G,M) \simeq \sha^2(\ell,\hat T_{E/\ell})$ and that $\sha^2_{\rm cycl}(G_{\rm prim},M_{\rm prim}) 
\simeq \sha^2(\ell,\hat  T_{E_{\rm prim}/\ell})$. 

\bigskip
\noindent
{\bf Proof of Theorem \ref {cyclic p}.}
By Theorem \ref {global}, we have 
$$\sha^1(\ell,T_{E/\ell}) = 0 \iff
 \sha^1(\ell,T_{E_{\rm prim}/\ell}) = 0,$$
 therefore, since $\hat T_{E/\ell} \simeq M$ and  $\hat T_{E_{\rm prim}/\ell} \simeq M_{\rm prim}$, we have 
 %$$\sha^2(\ell,\hat T_{E/\ell}) = 0 \iff
 %\sha^2(\ell, \hat T_{E_{\rm prim}/\ell}) = 0.$$
  $$\sha^2(\ell,M) = 0 \iff
 \sha^2(\ell, M_{\rm prim}) = 0.$$
 By Proposition \ref {lemma}, we have
 
 $$\sha^2_{cycl}(G,M)  \simeq  \sha^2(\ell,M)$$ and $$\sha^2_{cycl}(G_{\rm prim},M_{\rm prim})  \simeq  \sha^2(\ell,M_{\rm prim}).$$ Since
 $M = \hat T$ and $M_{\rm prim} = \hat T_{\rm prim}$, we obtain
 
 $$\sha^2_{\rm cycl}(G,\hat T) = 0 \iff \sha^2_{\rm cycl}(G_{\rm prim},\hat T_{\rm prim}) = 0,$$ as claimed.

\medskip
\noindent
{\bf Proof of Theorem \ref{torus p}.} Theorems \ref{cyclic p} and \ref{torus p} are equivalent by Theorem \ref {Colliot Br cycl}.

\bigskip
\noindent
From now on, we assume that the extensions $K_i/k$ are pairwise disjoint which implies the same property for the extensions $L_i/l$.

\bigskip
\noindent
{\bf Proof of Theorem \ref {prime cyclic sha}}
 Note that $G_{\rm prim}$ is an elementary abelian $p$-group, with $|G_{\rm prim}| =p$ if $n = 1$,  and $|G_{\rm prim}| \ge p^2$ if $n \ge 2$.
We may assume that $n \ge 2$.

\medskip
Assume first that $G_{\rm prim} \not \simeq C_p \times C_p$. Then  all the factors of the \'etale $\ell$-algebra $E$ are not contained in a field extension of degree $p^2$ of $\ell$,
therefore Theorem \ref{BLP} implies that $\sha^1(\ell,T) = 0$, and by Corollary \ref{0} this
implies that $\sha^2_{cycl}(G_{\rm prim},M) = 0.$

\medskip Assume now that $G_{\rm prim}  \simeq C_p \times C_p$.  Then, $\ell'$ is
a degree $p^2$ extension of $l$ containing all the factors of $E$.  Since the extensions $L_i$ are pairwise disjoint, $n \leq p+1$. Since $\ell'/\ell$ is unramified at all the finite places, the local degrees 
are $\le p$. By Theorem \ref{BLP}, this implies that $\sha^1(\ell,T) \simeq ({\bf Z} / p {\bf Z})^{n-2}.$ 
Hence $\sha^2(\ell,\hat T) ) \simeq ({\bf Z} / p {\bf Z})^{n-2}.$ 
Since the $G$-lattices $\hat T$ and $M$ are isomorphic, by  Corollary \ref{cycl F} we have $\sha^2_{\rm cycl}(G,M) \simeq  ({\bf Z} / p {\bf Z})^{n-2}.$

\medskip
\noindent
{\bf Proof of Theorem \ref {prime}.} Both (a) and (b) follow directly from Theorems \ref{Colliot Br cycl}, \ref{T = X} and \ref{prime cyclic sha}.

%the injection ${\rm Br}(X^c)/{\rm Im}({\rm Br}(k)) \to {\rm Br}(T^c)/{\rm Br}(k)$, cf. \S \ref {un}

\bigskip
\noindent
{\bf Proof of Theorem \ref {p+2 cycl}.} By Theorem \ref {prime cyclic sha}, we have $\sha^2_{\rm cycl}(G_{\rm prim},\hat T_{\rm prim})  = 0,$
and by Theorem \ref {cyclic p} this implies that $\sha^2_{\rm cycl}(G,\hat T) = 0$.

% \ref{BLP}, we have $ \sha^1(\ell,T_{E_{\rm prim}/\ell}) = 0$
%Theorem \ref {global}

\medskip
\noindent
{\bf Proof of Theorem \ref {p+2}.} Theorems \ref{p+2 cycl} and \ref{p+2} (a) are equivalent by Theorem \ref {Colliot Br cycl}.
Theorem \ref{p+2} (b) follows Theorem \ref{T = X}.
Theorem \ref{p+2} (c)  follows from   Corollary \ref{Sansuc} and Proposition \ref{Vosk}.

\bigskip
\noindent
{\bf Proof of Theorem \ref {cyclic at least one prime}.} By Theorem \ref {global at least one prime}, we have 
$$\sha^1(\ell,T_{E/\ell})^{\ast} \simeq
 \sha^1(\ell,T_{E_{\rm prim}/\ell})^{\ast},$$
 and by Poitou-Tate duality this implies that 
 $\sha^2(\ell, \hat T_{E/\ell}) \simeq
 \sha^2(\ell, \hat T_{E_{\rm prim}/\ell})$. Since $ \hat T_{E/\ell} \simeq M \simeq \hat T$, and 
 $\hat T_{E_{\rm prim}/\ell} \simeq M_{\rm prim} \simeq \hat T_{\rm prim}$, 
 applying  Corollary \ref{cycl F}, we obtain 
$\sha^2_{\rm cycl}(G_{\rm prim},\hat T_{\rm prim}) \simeq
\sha^2_{\rm cycl}(G,\hat T)$, as claimed.

% \ref{BLP}, we have $ \sha^1(\ell,T_{E_{\rm prim}/\ell}) = 0$
%Theorem \ref {global}

\medskip
\noindent
{\bf Proof of Theorem \ref {torus at least one prime}.} Theorems \ref{cyclic at least one prime} and \ref{torus at least one prime} are equivalent by Theorem \ref {Colliot Br cycl}.

\section{Unramified Brauer groups and products of cyclic extensions}\label{un cyclic}

\medskip Let $p$ be a prime number, and let 
%$k$ be a field of characteristic $\not = p$. 
$L = K_1 \times \dots \times K_n$, where $K_1, \dots, K_n$ are distinct cyclic extensions of $k$ of
degree $p$. Let $T_{L/{k}}  = R^{(1)}_{L/k}({\bf G}_m)$ be the $k$-torus
defined by $$1 \to T_{L/{k}} \to R_{L/k}({\bf G}_m) {\buildrel {\rm N}_{L/{k}} \over \longrightarrow}  {\bf G}_m.$$ 
Let $k'/k$ be a Galois extension of minimal degree splitting $T_{L/{k}}$, and let $G = {\rm Gal}(k'/k)$. Let $M$ be the $G$-lattice of characters of $T_{L/{k}}$.

\medskip
Let $a \in k^{\times}$, and let $X$ be the affine $k$-variety determined by the equation ${\rm N}_{L/k}(t) = a$, and note that $X$ is a
torsor under $T = T_{L/k}$. Let $X^c$ be
a smooth compactification of $X$.  

\medskip

The following result is an immediate consequence of Theorems  \ref{prime} and \ref{T = X} :

\begin{theo}\label{brun}

\medskip
{\rm (a)} If $G \not \simeq C_p \times C_p$, then $${\rm Br}(T^c)  / {\rm Br}(k) = {\rm Br}(X^c)  / {\rm Im}({\rm Br}(k) ) = 0.$$

\medskip
{\rm (b)} If $G  \simeq C_p \times C_p$, then $${\rm Br}(T^c)/{\rm Br}(k)  \simeq {\rm Br}(X^c)  / {\rm Im}({\rm Br}(k) ) \simeq ({\bf Z} / p {\bf Z})^{n-2}.$$

\end{theo}

%Assume now that {\it the characteristic of $k$ is not $p$.} In this case, we determine the structure of ${\rm Br}(X^c)  / {\rm Im}({\rm Br}(k) )$, as
%follows

%\begin{theo}\label {brun p} Suppose that ${\rm char}(k) \not = p$, and that $G  \simeq C_p \times C_p$. Then we have 
%$${\rm Br}(X^c)  / {\rm Im}({\rm Br}(k) ) \simeq ({\bf Z} / p {\bf Z})^{n-2}.$$

%\end{theo}

%\medskip
%The proof of this theorem is a combination of arithmetic and algebra. By Theorem \ref {brun} (b), we already know that 
%${\rm Br}(T^c)/{\rm Br}(k)  \simeq ({\bf Z} / p {\bf Z})^{n-2}.$ Recall from \S \ref{un} that ${\rm Br}(X^c)  / {\rm Im}({\rm Br}(k) )$ injects
%into ${\rm Br}(T^c)/{\rm Br}(k)$. On the other hand, Theorem \ref{generators} shows that ${\rm Br}(X^c)  / {\rm Im}({\rm Br}(k) )$
%is of dimension at least $n-2$ over $({\bf Z} / p {\bf Z})$.

\medskip Assuming that ${\rm char}(k) \not = p$,  we  obtain a more precise result, namely we give generators for the group  ${\rm Br}(X^c)  / {\rm Im}({\rm Br}(k))$. 

\medskip Set  $I = \{1, \dots, n\}$. 
We consider the norm polynomials ${\rm N}_{K_i/k}(t_i)$  for $i \in I$ as elements of $k(X)$.  For all $i \in I$, set $N_i = {\rm N}_{K_i/{k}}(t_i)$ and let $\sigma_i$ be a generator
of ${\rm Gal}(K_i/k)$. Let $\tilde K_n \in H^1(k,{\bf Z}/p{\bf Z})$ be the element associated to the pair $(K_n,\sigma_n)$. 
We identify $H^1(k(X^c),\mu_p)$ with $k(X^c)^*/k(X^c)^{*p}$ via the Kummer isomorphism and regard $[N_i] \in H^1(k(X^c),\mu_p)$.
The variety $X$ is defined by $N_1N_2 \dots N_n = a$ in the affine space $k[t_i, 1 \le i \le n]$. Let $(N_i,\tilde K_n)$ denote the class of the cyclic algebra of degree $p$ over $k(X^c)$ associated to $[N_i] \in H^1(k(X^c),\mu_p)$ and $\tilde K_n \in H^1(k(X^c), {\bf Z}/p{\bf Z})$. 

\medskip
\begin{theo}\label{generators} Suppose that ${\rm char}(k) \not = p$, and that $G  \simeq C_p \times C_p$. Then
 ${\rm Br}(X^c)  / {\rm Im}({\rm Br}(k))$  is generated by the  $n-2$ linearly independent elements $$(N_i,\tilde K_n), \ \ i = 1,\dots,n-2.$$ in ${\rm Br}(k(X^c))$. 

\end{theo}

\medskip

We begin with  the following lemma :

\begin{lemma}\label{kernel} Let $K/k$ be a cyclic extension of degree $n$ with $(n,{\rm char}(k)) = 1$. Let $\sigma$ be a generator of ${\rm Gal}(K/k)$
and let $A$ be the cyclic algebra over $k$ defined by $((K,\sigma),c)$ for some $c \in k^{\times}$. 
Let $X$ be the variety ${\rm N}_{K/k}(t) = c$. Then the kernel of ${\rm Br}(k) \to {\rm Br}(k(X))$ is generated by the class of $A$.
\end{lemma}

\noindent
{\bf Proof.} Let $Y_A$ be the Severi-Brauer variety of $A$. Since $A$ is split by $k(Y_A)$, the element $c$ is a norm from the extension
$K k(Y_A)/k(Y_A)$. Thus $X$ has a rational point over $k(Y_A)$  and the map ${\rm Br}(k(Y_A)) \to {\rm Br}(k(Y_A)(X))$ has trivial kernel. We
have, 
$${\rm Ker}({\rm Br}(k) \to {\rm Br}(k(X))) \subset {\rm Ker}({\rm Br}(k) \to {\rm Br}(k(X)(Y_A))) =$$
$$= {\rm Ker}({\rm Br}(k) \to {\rm Br}(k(Y_A))) = \langle [A] \rangle$$ by a theorem of Amitsur \cite{GS}, Theorem 5.4.1.
Since $[A]$ is zero in ${\rm Br}(k(X))$, it follows that
${\rm Ker}({\rm Br}(k) \to {\rm Br}(k(X)) = \langle [A] \rangle$ .

\bigskip
\noindent
In  the following proof  of Theorem \ref{generators},  we use the fact that the Brauer group of $X^c$ is the unramified Brauer group of $X^c$  (cf \cite {CT}, Section 2), namely the subgroup of ${\rm Br}(k(X^c))$  consisting of all elements which are unramified at all discrete valuations of $k(X^c)$.     This is a consequence of the purity results of Cesnavius \cite {C}, Theorem 1.2.

\bigskip
\noindent {\bf Proof of Theorem \ref{generators}.}  The strategy of the proof is the following. We first show that the algebras $(N_i,\tilde K_n), \ \ i = 1,\dots,n-2$ belong to $ {\rm Br}(X^c)$.
By Theorem \ref{brun} (b), we know that ${\rm Br}(T^c)/{\rm Br}(k)  \simeq ({\bf Z} / p {\bf Z})^{n-2}$. Moreover, ${\rm Br}(X^c)  / {\rm Im}({\rm Br}(k))$
injects into ${\rm Br}(T^c)/{\rm Br}(k)$ (see \S \ref{un}).
We next show that the elements
$(N_i,\tilde K_n), \ \ i = 1,\dots,n-2.$  are linearly independent over ${\bf Z}/p{\bf Z}$, and this yields the desired result.

\medskip
 Identifying $G$ with ${\bf Z}/p{\bf Z} \times {\bf Z}/p {\bf Z}$, let $\sigma(i,j) \in G$ correspond to $(i,j) \in {\bf Z}/p{\bf Z}
\times  {\bf Z}/p{\bf Z}$, $0 \le i \le p-1$, $0 \le j \le p-1$. We assume $K_1 = (k')^{\sigma(0,1)}$ and $K_n = (k')^{\sigma(1,0)}$. Pick a
generator $\sigma_1$ of ${\rm Gal}(K_1/k)$ and let $\tilde K_1 = [(K_1,\sigma_1)] \in H^1(k,{\bf Z}/p{\bf Z})$. For each $j \in I$, one can choose  a generator
$\sigma_j$ of ${\rm Gal}(K_j/k)$ such that   for $j \ge 2$, we have $\tilde K_j =  i_j \tilde K_1 + \tilde K_n$
for some $i_j$ with $1 \leq i_j \leq p-1$.
In fact for $j \geq 2$, we have $K_j = (k')^{\sigma(1,j)}$, and $\sigma_j$ is determined by $ \sigma_1$  and $\sigma_n$.

\bigskip
Since $N_i \in N_{K_i k(X^c)/k(X^c)}( K_i k(X^c))  $, we get $(N_i,\tilde K_i)_{k(X^c)} = 0$. Hence we have $(N_1,\tilde K_j) = (N_1,i_j\tilde K_1 + \tilde K_n) = (N_1,\tilde K_n)$
for every $j \in I$, $j \ge 2$.

We show that the cyclic algebras $(N_i,\tilde K_n)$ are unramified on $X^c$ for all $i \in I$, and
that $(N_1,\tilde K_n),\dots,(N_{n-2}, \tilde K_n)$ are linearly independent in ${\rm Br}(X^c)/{\rm Br}(k)$.

\medskip
Let $R$ be a discrete valuation ring containing the field $k$, with fraction field of $R$ equal to $k(X^c)$. We prove  that the algebras 
$(N_i,\tilde K_n)$ are unramified with respect to the valuation $v_R$. Let us denote by $\kappa$ the residue field of $R$, and let
$\partial_R : {\rm Br}(k(X^c)) \to H^1(\kappa,{\bf Q}/{\bf Z})$ be the residue map. 

\medskip
Let $\overline{[ K_i]}$ denote the image of $\tilde K_i$ in $H^1(\kappa,{\bf Z}/p{\bf Z}$). We have (see for instance \cite{GS} Lemma 6.8.4 and construction 6.8.5),

$$\partial_R(N_i,\tilde K_n) = \overline{[ K_n]}^{v_R(N_i)}.$$

%\medskipshekhar84112@gmail.com
%For all $i \in I$, we denote by  $[K_i] \in H^1(k(X^c),{\bf Z}/p{\bf Z})$ the class of the cyclic extension $k(X^c)K_i$ of $k(X^c)$ with 
%generator $\sigma_i$. Note that for all $i \in I$, we have either $[K_i ] = [K_1^jK_n]$ for some $j = 0,\dots, p-1$ or $[K_i] = [K_1]$.
%On the other hand, for all $i \in I$, the cyclic algebra $(N_i,K_i)$
%$(N_i,K_i)_{k(X^c)}$ 
%is trivial. Hence we have $(N_1,K_n) =
%(N_1,K_i)$ 
%$(N_1,K_n)_{k(X^c)} =(N_1,K_i)_{k(X^c)}$ 
%if $i \not = 1$. 

%\medskip
%Let  $\overline {[K_i]}$ be the image of $[K_i]$ in $H^1(\kappa,{\bf Z}/p{\bf Z})$. 
%We have $$\partial_R(N_i,K_n) = \overline {[K_n]}^{v_R(N_i)}.$$

\medskip
If $K_n \subset \kappa$  we have $\overline {[\tilde K_n]} = 0$ and $\partial_R(N_i,\tilde K_n) = 0$ for $1  \leq i \leq n-2.$.
Suppose that $\partial_R(N_1,\tilde K_n) \ne  0$.
Then $K_n$ is not contained in $\kappa$.  In this case, 
$ K_n \kappa$ is a degree $p$ cyclic extension of $\kappa$. 
%In the first case, we have $\partial_R(N_i,K_n) = 0$ for all $i \in I$. 
The extension $k(X^c)K_n/k(X^c)$ is cyclic of degree $p$, and has residual degree $p$, hence is unramified at $R$. Further
$N_n \in k(X^c)$ is a norm from the extension $k(X^c)K_n/k(X^c)$.
%$N_n$ is an element of $k(X^c)$ and it is a norm from the corresponding extension. 
Hence the valuation $v_R(N_n)$ is divisible by $p$. Since $(N_1, \tilde K_j) = (N_1,\tilde K_n)$ for $j \ge 2$, we have 
$\partial_R(N_1,\tilde K_j) \not = 0$, and $K_j$ is not contained in $\kappa$. 
Repeating the above argument,  we see that $p$ divides $v_R(N_j)$ for all $j \ge 2$. Since $a = N_1 \dots N_n \in k$, $v_R(a) = 0$, 
 $p$ divides $v_R(N_j)$ for $2 \le j \le n$, and it follows that $p$ divides $v_R(N_1)$. This implies that 
$\partial_R(N_1,\tilde K_n) = (\overline K_n)^{v_R(N_1)}= 0$, contradicting the assumption that $\partial_R(N_1,\tilde K_n)  \not = 0$. 
This implies that $\partial_R(N_1,\tilde K_n)  = 0$. A similar argument, interchanging $1$ and $i$, with $i \le n-1$, gives that
$\partial_R(N_i,\tilde K_n)  = 0$. Hence the elements $(N_i,\tilde K_n)$ are unramified at $R$ for every discrete valuation ring $R$
with field of fractions $k(X^c)$. By purity for ${\rm Br}(X^c)$, we have 
$(N_i,\tilde K_n) \in {\rm Br}(X^c)$. 

%with
%$K_n$ replaced by $K_i$ for $i \not = 1$, we see that $p$ divides $v_R(N_i)$ for all $i = 2,\dots,n$. 

%\medskip
%We have $a = N_1 \dots N_n$, hence $p$ also divides $v_R(N_1)$. We have
 %$$\partial_R(N_1,K_n) = \overline {[K_n]}^{v_R(N_1)} = 0,$$ since $H^1(\kappa,{\bf Z}/p{\bf Z})$ is a $p$-torsion group. Therefore
 %the elements $(N_i,K_n)$ are unramified with respect to $R$ for all $i \in I$. 

%\medskip
%Let $E$ be the composite field of the subfields $K_1,\dots,K_n$; by hypothesis, $E/k$ is Galois and ${\rm Gal}(E/k)$ is isomorphic
%to $C_p \times C_p$. 
%The classes of the degree $p$ subfields of $f(X^c)E$ in $H^1(k(X^c),{\bf Z}/p{\bf Z})$ are $[K_1]$ and $[K_1]^j[K_n]$ for
%$j = 0,\dots,p-1$. If $K_i$ is a degree $p$ subfield of $E$, then $[K_i] = [K_1]^j[K_n]$ for some $j = 0,\dots,p-1$ or $[K_1]$. Note that
%$(N_1,K_i) = (N_1,K_1^jK_n) = (N_1,K_n)$; in this case, either the residue is 0, or $p$ divides $v_R(N_i)$. 

 \medskip
 Let us check that the algebras $(N_1,\tilde K_n),\dots,(N_{n-2}, \tilde K_n)$ are linearly independent  in ${\rm Br}(k(X^c))/{\rm Br}(k)$. Let us
 project $X$ to the $d$-dimensional affine space, where $d = (n-1)p$, corresponding to the coordinates involving the first $n-1$ norm polynomials. Let 
 $M$ be the function field of this affine space; we have $k \subset M \subset k(X^c)$. Note that $N_1,\dots,N_{n-1} \in M$. 
 We have $(N_i, \tilde K_n) \in {\rm Br}(M)$ for all $i = 1,\dots,n-1$.

 \medskip We want to show that the algebras $(N_1,\tilde K_n),\dots,(N_{n-2},\tilde K_n)$ are linearly independent in ${\rm Br}(k(X^c))/{\rm Br}(k)$. If not, then there
 exist $r_1,\dots,r_{n-2} \in {\bf Z}$ not all zero 
 with $0 \le r_i \le p-1$ such that $\underset{i = 1,\dots,n-2}\sum r_i (N_i,\tilde K_n)_{k(X^c)} = \alpha$ for some $\alpha \in {\rm Br}(k)$. 
 
 \medskip The kernel of the natural homomorphism ${\rm Br}(M) \to {\rm Br}(k(X^c))$ is generated by the class of the algebra $(a^{-1}N_1\dots N_{n-1},\tilde K_n)$, 
 by  Lemma \ref{kernel}   applied to $c = aN_1^{-1}. N_2^{-1}.....N_{n-1}^{-1}$.

 Hence there exists $s \in {\bf Z}$ with $0 \le s \le p-1$ such that 
 
 $$\underset{i = 1,\dots,n-2}\sum r_i (N_i,\tilde K_n)_{M} - \alpha = s(a^{-1}N_1\dots N_{n-1},\tilde K_n)_{M}$$ in ${\rm Br}(M).$ The polynomials $N_i$ are irreducible (see \cite {Fla}, Theorem 2). 
  Take residue on both sides at the valuation $v_{N_i}$ corresponding
 to the irreducible polynomial $N_i$ for all $i = 1, \dots, n-2$. The residue of the left side is $r_i \overline{[K_n]}$,
 and the right side is $s \overline{[K_n]}$.
 
 \medskip
 {\bf Claim.} $ \overline{[K_n]} \not = 0$ in the residue field $\kappa(v_{N_i})$ for all $i = 1,\dots,n-2$. 
 
 \medskip Assume that the claim holds. Then we have $r_i = s \neq 0$ for all $i = 1,\dots,n-2$, and we get $s(a^{-1}N_{n-1},\tilde K_n) = - \alpha$. With
 respect to the valuation $v_{N_{n-1}}$, we
 have the residue $\partial (a^{-1}N_{n-1},\tilde K_n) =   \overline{[K_n]} \not = 0$ by the claim, and this leads to a contradiction since $\partial(\alpha)=0$, $\alpha$ being an element of $k$,  and $s \neq 0$.
 
 \medskip
 It remains to prove the claim. Let us show that $K_n$ is not contained in  $\kappa(v_{N_i})$ for all $i = 1,\dots,n-2$. Let $M_i$ be the
 function field of the $k$-variety determined by the polynomial $N_i$.
 %and let us show that $K_n$ is not contained in $\kappa(v_{N_i})$. 
 Since $\kappa(v_{N_i})$ is a purely transcendental extension of $M_i$ obtained by adjoining the coordinates involved in the polynomials  $\{N_j, j \neq i, j \leq n-1\} $,  it suffices to show that $K_n$ is not contained in $M_i$. Suppose
 that $K_n$ is a subfield of $M_i$. Let us base change to $K_i$ : 
 %we have $E = K_iK_n$, 
 the field $K_iK_n$ is a subfield of $K_iM_i$. 
 After base change to $K_i$, the polynomial $N_i$ is transformed to the product  $X_1\dots X_p$ for some variables $X_1,\dots,X_p$. 
 Hence $K_iM_i$ is a product of rational function fields over $K_i$; therefore it cannot contain $K_n$, thereby leading to contradiction.

 \section {Some consequences for semi-global fields} \label{semiglobal}
 
 Let $K$ be a complete discrete valued field with 
valuation ring $O$ and residue field $\kappa$. Let $X/K$ be a normal projective, geometically integral 
curve over $K$ and let $F = K(X)$. We call $F$ a {\it semi-global field}. In \cite {CTPS}, \S 2.3, certain higher reciprocity obstructions 
were constructed to study  the failure of the Hasse principle for varieties over $F$ with respect to discrete valuations
of $F$ centered on a regular proper model $\mathcal X / O$ of the curve $X/K$. It was also proved in \cite {CTPS}, Example 2.6, that  these obstructions do not suffice to detect the failure of the Hasse principle on principal homogeneous spaces under tori  defined over $F$. Using an example of a multinorm  torus constructed in \cite{Sumit}, Corollary 7.12, we give another example where  the  obstructions constructed
in \cite {CTPS} do not suffice to detect failure of Hasse principle.

\medskip
We recall the following construction from \cite{Sumit}, Corollary 7.12. 
Let $K = {\bf C}((t))$ and $F = {\bf C}((t))(X)[Y]/(XY(X+Y-1)(X-2)-t)$. Let $L_1 = F((XY)^{1/n}, (Y(X + Y - 1))^{1/n})$ and $L_2
= F((XY \theta_1)^{1/n},(Y(X + Y - 1)\theta_2)^{1/n})$, where $\theta_1 = (X-2)/(X-2+XY(X+Y-1))$ and $\theta_2 = (Y-2)/(Y-2+XY(X + Y -1)$. Then $L_1$ and
$L_2$ are Galois extensions of $F$ which are linearly disjoint over $F$ (c.f \cite{Sumit}, Corollary 7.12).
%constructed in \cite{CTPS}, and let 
Let $T = {\rm R}^1_{L_1\times L_2/F}(\bf G_m)$ be the associated norm one torus. It is proved in
 \cite{Sumit}, Corollary 7.12,  there is a principal homogeneous space under $T$
which fails the Hasse principle with respect to all discrete valuations of $F$.  The proof invokes $R$-equivalence of tori to prove
that Hasse principle fails in the patching setting of Harbater-Hartmann-Krashen. The second step is to prove that
failures of Hasse principle in the patching setting implies the failure of the Hasse principle with respect to all discrete valuations
of $F$. Since the cohomological dimension of $F$ is 2, 
the only  obstruction of \cite {CTPS} in this case is 
the one coming from the Brauer group of $X^c$. In view of  Theorem \ref{Polliot-Rapinchuk}, we have ${\rm Br}(X^c)/{\rm Br}(k) = 0$, and the obstruction vanishes.
However, Hasse principle fails for $X$.


\begin{thebibliography}{99}
\bibitem[BLP 19]{BLP} E. Bayer-Fluckiger, T-Y. Lee, R. Parimala, \textit{Hasse principles for multinorm equations}, Adv. Math. \textbf {356} (2019), 35 pages.
\bibitem[Bo 91]{Bo}  A. Borel, \textit{Linear algebraic groups} Second edition. Graduate Texts in Mathematics, {\bf 126}, Springer-Verlag, New York (1991).
\bibitem[C 19]{C} K. Cesnavicius, \textit{Purity for the Brauer group}, Duke Math. J. {\bf 168} (2019), 1461-1486.
\bibitem[CT 14]{CT} J-L. Colliot-Th\'el\`ene, \textit{Groupe de Brauer non ramifi\'e d'espaces homog\`enes de tores}, J. Th\'eor. Nombres  Bordeaux {\bf 26} (2014), 69-83.
\bibitem[CTHSk 03]{CTHSk 03} J-L. Colliot-Th\'el\`ene, D. Harari,  A. Skorobogatov, \textit{Valeurs d'un polyn\^ome \`a une variable repr\'esent\'ees par une norme},  Number theory and algebraic geometry,  London Math.Soc.Lecture Note Ser {\bf 303},Cambridge Univ.Press,Cambridge (2003), 69-89.
\bibitem[CTHSk 05]{CTHSk 05} J-L. Colliot-Th\'el\`ene, D. Harari,  A. Skorobogatov, \textit{Compactification \'equivariante d'un tore} (d'apr\`es Brylinski
et K\"unnemann), Expo. Math {\bf 23} (2005), 161-170.
\bibitem[CTPS 16]{CTPS} J-L. Colliot-Th\'el\`ene, R. Parimala, V. Suresh, \textit{Lois de r\'eciprocit\'e sup\'erieures et points rationnels} (French) [Higher reciprocity laws and rational points] Trans. Amer. Math. Soc. {\bf 368}  (2016), 4219-4255. 
\bibitem[CTSk 19]{CTSk} J-L. Colliot-Th\'el\`ene, A. Skorobogatov, \textit {The Brauer-Grothendieck group} (2019), available on the authors' web page.
\bibitem[CTS 77]{CTS 77} J-L. Colliot-Th\'el\`ene, J-J. Sansuc, \textit{La R-\'equivalence sur les tores}, Ann. Sc ENS \textbf {10} (1977),
175-229.
\bibitem[CTS 87]{CTS 87} J-L. Colliot-Th\'el\`ene, J-J. Sansuc, \textit{Principal homogeneous spaces under flasque tori~: applications}, J. Algebra \textbf {106} (1987),
148-205.
\bibitem[F 53]{Fla} H. Flanders, \textit{The norm function of an akgebraic field extension}, Pacific J. Math. \textbf{3} (1953), 103-113.
\bibitem[F 62]{F} A. Fr\"ohlich, \textit{On non-ramified extensions with prescribed Galois group}, Mathematika \textbf{9} (1962), 133-134.
\bibitem[GS 06]{GS} P. Gille, T. Szamuely,
\textit{Central simple algebras and Galois cohomology}, Cambridge University Press (2006).
\bibitem[HS 71]{HS}, P. Hilton, U. Stammbach, \textit{A course in homological algebra}, Graduate Texts in Mathematics,
\textbf {4}, Springer Verlag (1971).
\bibitem[H 84]{H} W. H\"{u}rlimann,
\textit{On algebraic tori of norm type}, Comment. Math. Helv.  \textbf{59} (1984), 539-549.
\bibitem[Ma 19]  {Ma} A. Macedo, On the obstruction to the Hasse principle for multinorm equations,  arXiv:1912.11941v1.
\bibitem[PR 13]{PR} T. Pollio, A. Rapinchuk, \textit{The multinorm principle for linearly disjoint Galois extensions}, J. Number Th. \textbf {133}, (2013), 802-821.
\bibitem[San 81]{San} J.-J. Sansuc, \textit{Groupe de Brauer et arithm\'etique des groupes alg\'ebriques lin\'eaires sur un corps de nombres},
J. Reine Angew. Math. \textbf{327} (1981), 12-80.
\bibitem[S 79]{Se} J-P. Serre, \textit {Local fields}, Graduate Texts in Mathematics {\bf 67}, Springer Verlag (1979).
\bibitem[Su 19]{Sumit} C.M. Sumit, \textit{Local-global principle for norm one tori}, arXiv:1904.00966v1.
\bibitem[V ]{V} V. E. Voskresenskii, \textit{Algebraic groups and their birational invariants}, Trans. of Math. Monographs \textbf {179} AMS (1998).

\end{thebibliography}
\end{document}